\documentclass[11pt,a4paper]{article}

\usepackage{setspace}
\usepackage[utf8]{inputenc}
\usepackage[T1]{fontenc}
\usepackage[english]{babel}

\usepackage{amsmath, amssymb, amsthm}
\allowdisplaybreaks 
\usepackage{mathtools}
\usepackage{geometry}
\usepackage{booktabs}
\usepackage{xcolor}

\usepackage{physics}
\usepackage{mathtools}
 \usepackage{algorithm}
\usepackage{algpseudocode}

\usepackage{xcolor}

\usepackage{sidecap}

\usepackage{setspace}  

\newcounter{rupertcommentno}
\definecolor{piniengruen}{RGB}{46,139,87}

\newcounter{julianecommentno}
\usepackage{graphicx}
\usepackage{subfigure}
\usepackage{caption}
\usepackage{hyperref}
\hypersetup{
    colorlinks=true,
    linkcolor=black,
    citecolor=black,
    filecolor=black,
    urlcolor=black
}
\usepackage{enumitem}
\usepackage{cleveref}
\crefname{figure}{Figure}{Figures}

\usepackage{lineno} 
\usepackage{physics}

\newtheorem{theorem}{Theorem}[section]

\theoremstyle{remark}
\newtheorem{remark}[theorem]{Remark}

\title{A WKB-related time-stepping scheme for differential equations describing  oscillatory systems}
\author{Juliane Rosemeier \\ \small Freie Universit\"at Berlin \\ \small \texttt{juliane.rosemeier@fu-berlin.de}
\and 
Rupert Klein \\ 
\small Freie Universit\"at Berlin \\ 
\small \texttt{rupert.klein@math.fu-berlin.de}
    }
\date{\today}

\begin{document}

\maketitle

\begin{abstract}
In this study, we present a novel time-stepping scheme for multiscale  differential equations describing oscillatory systems with well-separated scales, where the scale separation is controlled by a small parameter $\epsilon$. The time-stepping method is related to a multi-modal WKB approximation and relies on a transformation of variables derived in this work.  The analysis reveals that, in the transformed formulation, the leading-order oscillations are either eliminated or appear only at higher asymptotic order.
The method is applied to ordinary differential equations, including the well-known van der Pol oscillator.
We investigate the accuracy of the proposed method numerically for different parameter regimes, in particular for decreasing values of $\epsilon$, and study how the parameters of the numerical scheme must be adapted as $\epsilon$ is reduced. In the presented numerical tests, the computational cost remains bounded as $\epsilon$ is decreased.

\end{abstract}


\section{Introduction}

In this study, we propose  novel numerical time-stepping methods related to a \emph{WKB (Wentzel-Kramers-Brillouin) ansatz} to solve differential equations governed by fast oscillations, that is, problems that exhibit oscillatory stiffness. Such differential equations arise in a wide range of applications, including
geophysical and astrophysical fluid dynamics, \cite{Majda2003,Klein2010}, classical and quantum mechanics, \cite{LasserLubich2020,Nolting2014},
plasma physics, \cite{BirdsallLangdon1991,Nicholson1983}, and electromagnetic wave propagation, \cite{Griffiths2017,Jackson2013}. The highly oscillatory nature of these nonlinear equations presents significant challenges for standard numerical methods. Faster oscillations require smaller time steps to ensure stability of the solutions and achieve a given accuracy tolerance. 

We develop the numerical time-stepping schemes in the context of ordinary differential equations.
Ordinary differential equations describing oscillatory systems serve as a natural starting point for the study of hyperbolic and wave-type PDEs, since their discretizations frequently lead to coupled ODE systems that arise from oscillatory models. These problems are challenging both analytically and numerically due to the presence of fast temporal scales and predominantly imaginary spectra.


Numerical time integration of  differential equations describing oscillatory systems has attracted significant attention in recent years. 
A variety of approaches have been proposed to efficiently resolve the fast oscillations and to study the underlying dynamical systems. 
These include  the  multiscale methods in \cite{VaterKlein2018, ArtsteinLinshizTiti2007, ArtsteinGearKevrekidisSlemrodTiti2011}, integrating factor and exponential integrator methods \cite{HochbruckLubich1997,CoxMatthews2002}, which treat stiff linear oscillatory terms analytically, 
splitting methods such as Strang splitting \cite{Strang1968}, modulation-based techniques like modulated Fourier expansions \cite{HairerLubich2000,HairerLubichWanner2006}, 
multiscale approaches such as the heterogeneous multiscale method (HMM) \cite{EEngquist2003,SharpTsaiEngquist2005}, and parallel-in-time algorithms such as Parareal \cite{LionsMadayTurinici2001}, MGRIT \cite{Falgout2014MGRIT}, and PFASST \cite{EmmettMinion2012}. The study \cite{Gander2015} provides an overview of parallelizable time-stepping schemes. Time-parallel methods, like Parareal, frequently struggle when applied to oscillatory and hyperbolic problems, see for instance \cite{GanderVandewalle2007,GanderWuZhou2025}, but there is recent progress, see for example \cite{KrzysikDeSterckFalgoutSchroder2025,HowseDeSterckFalgoutMacLachlanSchroder,HautWingate2014,PeddleHautWingate2019,RosemeierHautWingate2024}.  
 In general, numerical methods for the highly oscillatory problems pursue different strategies such as integrating the fast linear part exactly, using micro-simulations, or representing the solution in a Fourier basis. 
For the method proposed in this study, micro-simulations are performed over short time windows to extract amplitude–phase coordinates from the solution. These data are then extrapolated to advance the solution over larger time steps, after which new micro-simulations are carried out.
A key novelty of our approach lies in the use of an amplitude–phase coordinate transformation with a nonlinear phase function.
While nonlinear phase functions have also been considered in \cite{KleinPeters1988, LasserLubich2020}, many existing approaches rely on linearized phase approximations.

For the application of the proposed numerical methods, the right-hand side of the differential equation is not needed in the representation of the new coordinates. The entire methodology is composed of a standard numerical time-stepping schemes applied to the problem in Cartesian coordinates in  short time windows and data analysis. In the present work, however, we also derive the governing equations in the transformed coordinates and analyze the resulting system, as this provides additional insight into the underlying dynamics.

\begin{figure}
    \centering
    \includegraphics[width=0.8\linewidth]{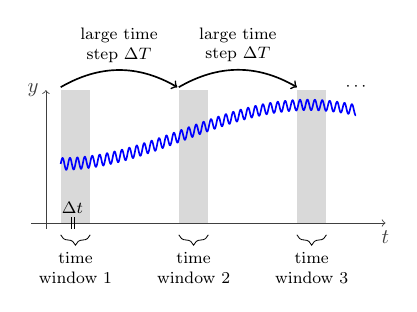}
    \caption{Illustration of a time-stepping method using a WKB ansatz. The blue line represents the solution of a differential equation describing an oscillatory system.
    A numerical approximation is sought for this problem. The idea of the  WKB-related time-stepping method is as follows: We compute solutions using a standard time-stepping scheme in short time windows with a small time step $\Delta t$. The data is projected onto amplitude–phase coordinates, and data-driven estimates of their slow evolution are constructed. Then,  the WKB approximation is used to take a large time step $\Delta T$ without applying a standard time-stepping scheme between the time windows. This has the potential to reduce the computational effort for time-stepping, although it requires additional operations.}
    \label{fig:idea_WKB_method}
\end{figure}

A drawback of standard time-stepping schemes, \textit{e.g.} Runge-Kutta methods, when applied to problems with oscillatory stiffness is that small time steps must be used throughout the entire solution domain to resolve the rapid oscillations, making computations expensive. 
Motivated by this observation, we aim to exploit analytical knowledge of the oscillatory behavior to develop new time-stepping schemes that reduce the computational cost. 
By performing numerical integration over a short time interval and leveraging the oscillatory structure of the solutions, we can advance the solution to a much later time beyond this short interval, \textit{i.\ e.}, we take a large time step.
By iterating this process, we obtain an approximation over the entire solution domain (see Figure \ref{fig:idea_WKB_method}). 

The new time-stepping schemes are related to WKB approximations.
A WKB ansatz is a common technique in physics and is well established in
theoretical fluid dynamics; see, for example, \cite{sanders_etal_2007, Majda2003, AchatzEtAl2010, SchlutowEtAl2017}. It can be used to construct approximations to differential equations and is
particularly suited for oscillatory problems. In this approach, the solution
is represented as a superposition of rapidly oscillating wave modes with a
potentially nonlinear phase and a slowly varying amplitude. The ansatz that we will use is developed in detail in \cref{subsec:general_Weakly_nonlinear_system,sec:Example_systems,sec:WKB-related_time-stepping scheme}. Since the functions describing the phases and amplitudes in the problems under investigation are, to leading order, slow, they are well suited for numerical integration with large time steps, which will be exploited in the construction of the WKB-based time-stepping schemes.



We consider two systems. In the first system, the fast oscillation completely vanish through the transformation. In the second system, they are still present, but appear at high asymptotic order, \textit{i.\ e.} not at the leading asymptotic order. This has important consequences for the numerical schemes proposed.

The remainder of the paper is organized as follows. \Cref{subsec:general_Weakly_nonlinear_system} introduces the model problem and derives the corresponding transformation. In \cref{sec:Example_systems}, we discuss two ordinary differential equations that fit into this framework and demonstrate their analysis. In \cref{sec:WKB-related_time-stepping scheme}, we introduce the new time-stepping scheme. Numerical experiments are presented in \cref{sec:Numerical_experiments}. Finally, we conclude with a discussion of the results and an outlook in \cref{sec:Discussion_Outlook}.


\section{Weakly nonlinear system and derivation of the transformation}
\label{subsec:general_Weakly_nonlinear_system}

The general form of the weakly nonlinear system under investigation is given by
\begin{equation}
\label{eq:general_problem}
     \vb{\dot y} = \frac{1}{\epsilon} \mathcal{L} \vb{y} + \mathcal{N}(\vb{y}) ,
\end{equation}
where $\epsilon > 0$ is a small parameter, $\mathcal{L}$ is a skew-Hermitian linear operator and $\mathcal{N}$ is a nonlinear term. In this exposition, we consider ordinary differential equations and, in particular, linear operators represented by real-valued matrices. More specifically, we assume that $\mathcal{L}$
 is real and skew-symmetric. In future work, more complex systems, in particular partial differential equations, will be considered, where the operators are not restricted to real-valued entries.

The corresponding linear problem is defined by
\begin{equation}
\label{eq:linaer_problem}
     \vb{\dot y}_{\text{lin}} = \frac{1}{\epsilon} \mathcal{L} \vb{y}_{\text{lin}} .
\end{equation}

For the linear problem we can compute the exact solution. For real-valued initial conditions, the solution remains real-valued.  Assuming that $\mathcal{L}$ is a real-valued matrix we find
\begin{equation}
\label{eq:lin_sol}
    \vb{y}_{\text{lin}} = c_0 \vb{r}_{0} + \sum_{j=1}^J c_{j} e^{\frac{i}{\epsilon}\mu_j t} \vb{r}_{j} + \overline{c}_{j} e^{-\frac{i}{\epsilon} \mu_j t} \overline{\vb{r}_{j}} ,
\end{equation}
where $i$ denotes the imaginary unit, the set $\sigma \backslash \{0\} = \{ \pm i \mu_j | j = 1, \dots, J\} \backslash \{0\}$ is the spectrum of $\mathcal{L}$ without 0, $\vb{r}_{j}, \overline{\vb{r}_{j}}$ denote the corresponding eigenvectors, and $c_{j}, \overline{c}_{j}$ are the coefficients determined by the initial conditions. If 0 is in the spectrum of the linear operator $\mathcal{L}$, the vector $\vb{r}_{0}$ is the corresponding eigenvector.
If $0$ is not an eigenvalue of $\mathcal{L}$, the first term in the sum on the right-hand side of \eqref{eq:lin_sol} is omitted.
This presentation of the solution ensures that it remains real valued. (Note: for a real-valued skew-symmetric matrix $\mathcal{L}$ the eigenvectors to complex-conjugate eigenvalues are complex-conjugate, too.)
 Writing $c_j$ in polar form, we obtain $c_j = a_j e^{i \psi_j}$ and equation \eqref{eq:lin_sol} becomes
 \begin{equation}
 \label{eq:genaral_lin_sol}
    \vb{y}_{\text{lin}} = a_0 \vb{r}_{0} + \sum_{j=1}^J  \left( a_j e^{i \psi_j}  e^{\frac{i}{\epsilon}\mu_j t} \vb{r}_{j} + a_j e^{-i \psi_j} e^{-\frac{i}{\epsilon} \mu_j t} \overline{\vb{r}_{j}} \right) .
\end{equation}

To account for the nonlinear term in \eqref{eq:general_problem}, we allow
$a_j$ and $\psi_j$ to depend on time. This step is also pursued in the \emph{method of variation of parameters}; see, for instance, \cite{sanders_etal_2007} (p.~16), where it leads to a standard form representing
a certain problem in new coordinates, see appendix \ref{app:derivation_standard form}.
This procedure leads to the following representation of the problem

\begin{equation}
\label{eq:general_WKB_ansatz}
    \vb{y} = a_0(t) \vb{r}_{0} + \sum_{j=1}^J  a_j(t)  e^{\frac{i}{\epsilon}( \mu_j t+\epsilon\psi_j(t))} \vb{r}_{j} + a_j(t) e^{-\frac{i}{\epsilon}(\mu_j t+\epsilon\psi_j(t))} \overline{\vb{r}_{j}} .
\end{equation}
where  $a_j(t) $ and $\psi_j(t)$ are real-valued functions depending on $t$. For simplicity, we will also employ the definition
\begin{equation}
\label{eq:def_phi}
\phi_j(t) = t+\epsilon \psi_j(t) .
\end{equation}
The representation \eqref{eq:general_WKB_ansatz} with the coordinates $a_j(t)$ and $\psi_j(t)$ can also be interpreted as a type of standard form. A crucial difference to the standard form in appendix \ref{app:derivation_standard form} is that the method of variation of parameters is applied to parameters expressed in polar rather than Cartesian coordinates.


The eigenvectors of the linear operator $\mathcal{L}$ form a basis of $\mathbb{C}^n$, 
\textit{i.\ e.}, $\vb{r}_j$ and $\overline{\vb{r}_j}$ (and $\vb{r}_0$ if $0$ is an eigenvalue) 
in equation \eqref{eq:WKB_ansatz_van_der_Pol}.
The coefficients of the eigenvectors  are expressed with the functions $a_j(t)$ and $\phi_j(t)$, which are chosen such that the components of  $\vb{y}(t)$ are real-valued. In particular, $a_j(t)$ and $\frac{\phi_j(t)}{\epsilon}$ are the polar coordinates of the coefficients. Consequently, with appropriate choices for  $a_j(t)$ and $\phi_j (t)$,  the solution of  \eqref{eq:general_problem}  can be described exactly. That is, $a_j(t)$ and $\phi_j (t)$ define a transformation.

Note that for a complex-conjugate eigenvalue pair $\pm i \mu_j$, if the vector $\vb{y}$ has no component along the corresponding eigenvectors $\vb{r}_j$ and $\overline{\vb{r}_j}$, the amplitude $a_j$ vanishes and the phase $\phi_j$ is no longer uniquely defined.


\section{Example systems}
\label{sec:Example_systems}

In this section, we consider two ordinary differential equations of the form \eqref{eq:general_problem}.  
For the first problem, the linear operator has the eigenvalue $0$, and after applying the transformation described in the previous section, the problem admits a form that is very favorable for numerical simulation. We denote this system as the problem with $0$-eigenvalue. The second problem is the van der Pol oscillator. Its structure is not as simple as that of the first problem after the transformation, which has important consequences for the time-stepping method we consider later.

\subsection{System 1: problem with $0$-eigenvalue}

\begin{figure}[h]
    \centering
    \includegraphics[width=1.0\linewidth]{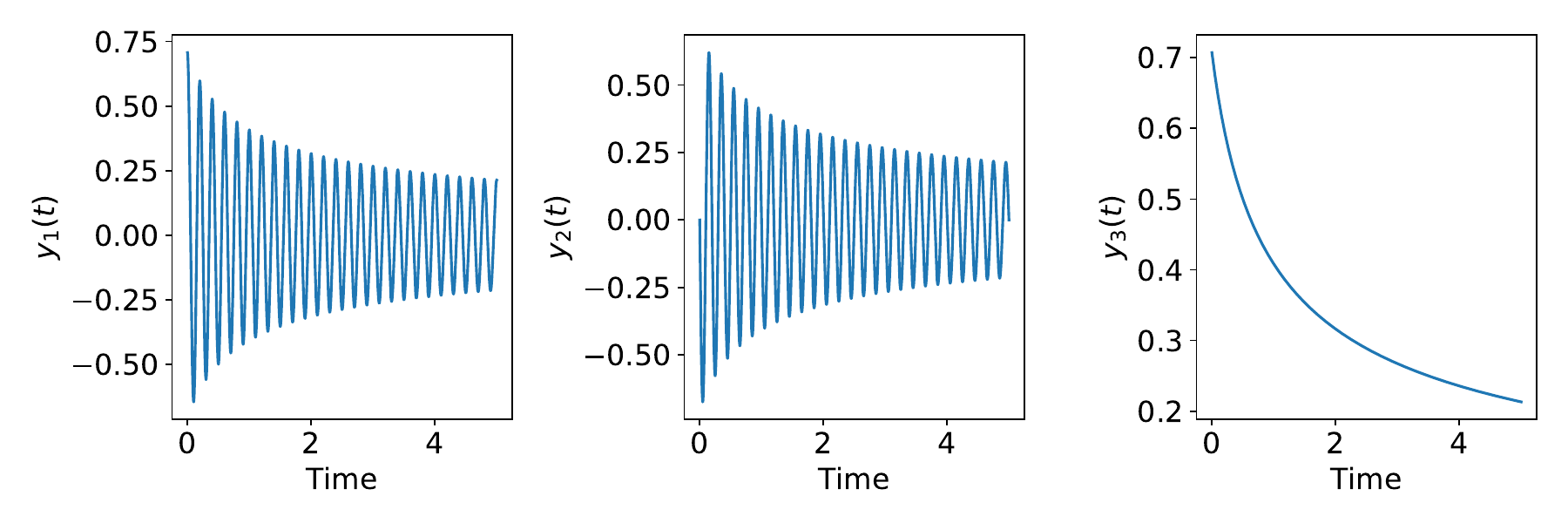}
    \caption{Problem with $0$-eigenvalue. The initial value is given by $\vb{y}_0 = (\sqrt{2}/2, 0.0, \sqrt{2}/2)^T$ and the parameters are chosen as $\epsilon = \frac{2.0 \times 10^{-1}}{2 \pi} \text{ and } \lambda = 1.0$. The solution in  Cartesian coordinates $y_1, y_2$ and $y_3$ is displayed.}
    \label{fig:problem_with_0_eigenvalue_cartesian}
\end{figure}

As a first test case, we consider the following system of ordinary differential equations:

\begin{equation}
\label{eq:ODE_with_zero_frequency}
     \vb{\dot y} = -\frac{1}{\epsilon} \begin{pmatrix}
        0 \\0 \\ 1
    \end{pmatrix} \times \vb{y} -
    \lambda \|\vb{y}\|^2 \vb{y} , \qquad \lambda > 0, \ \lambda \in O(1) ,
\end{equation}
for $\vb y = (y_1,y_2,y_3)^T \in \mathbb{R}^3$.
The system admits the form of  equation \eqref{eq:general_problem} with
\begin{equation}
\label{eq:ODE_with_zero_frequency_2}
    \mathcal{L} = \begin{pmatrix}
        0 & 1 & 0 \\
        -1 & 0 & 0 \\
        0 & 0 & 0 
    \end{pmatrix}
    , \qquad \mathcal{N} (\vb{y}) = -\lambda (y_1^2 + y_2^2 + y_3^2) \begin{pmatrix}
        y_1 \\ y_2 \\ y_3
    \end{pmatrix}
\end{equation}
The eigenvalues are $\mu_0 = 0$, $\mu_{1} =  i$ and $\overline{\mu}_{1} =  -i$ , with the corresponding eigenvectors $\vb{r}_0 = (0,0,1)^T$, $\vb{r}_{1} = \frac{1}{\sqrt{2} }
( 1, - i, 0 )^T$ and $\overline{\vb{r}}_{1} = \frac{1}{\sqrt{2} }
( 1, i, 0 )^T$ . \Cref{fig:problem_with_0_eigenvalue_cartesian} shows a solution of system \eqref{eq:ODE_with_zero_frequency}. In particular, we observe rapid oscillations in the $y_1$- and $y_2$-components.

The linear part of this equation has a zero eigenvalue. Such spectral features appear in various oscillatory systems with slow modes. As an example this can be motivated by structures found in geophysical fluid dynamics, such as Rossby waves.  The nonlinear term mimics weakly nonlinear interactions, making this test case relevant for a broad class of problems.

\begin{figure}[h]
    \centering
    \includegraphics[width=1.0\linewidth]{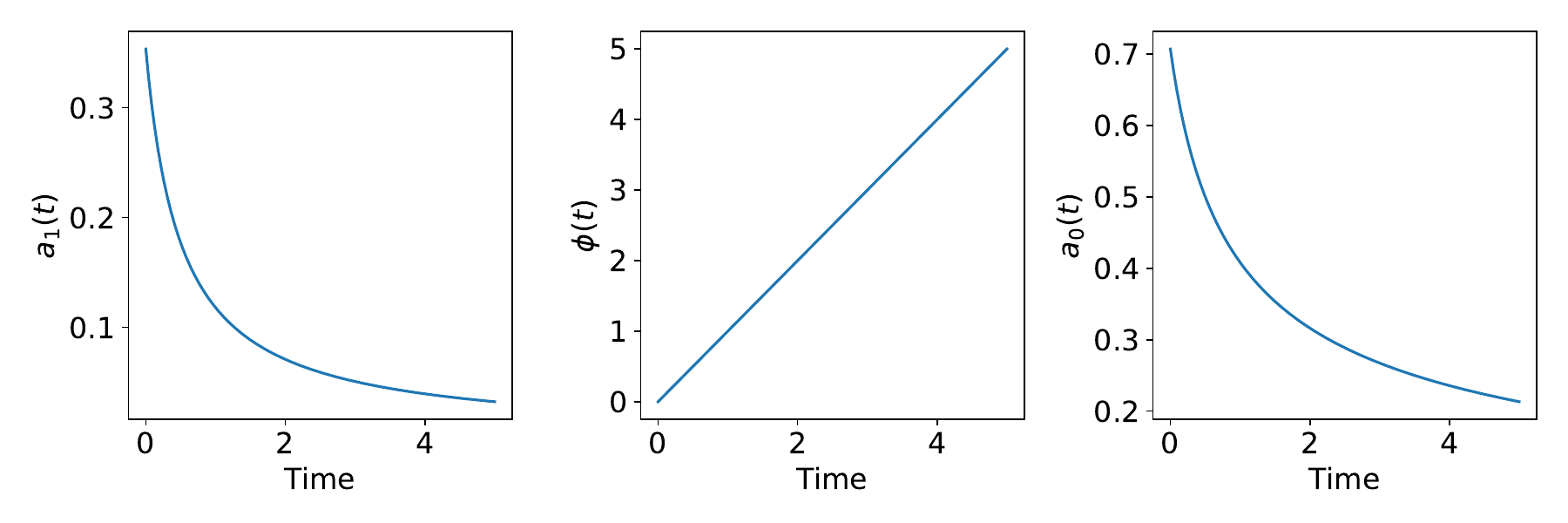}
    \caption{Problem with $0$-eigenvalue. The initial value is given by $\vb{y}_0 = (\sqrt{2}/2, 0.0, \sqrt{2}/2)^T$ and the parameters are chosen as $\epsilon = \frac{2.0 \times 10^{-1}}{2 \pi} \text{ and } \lambda = 1.0$. The solution is shown in the new  coordinates $a_0, a_1$ and $\phi_1$.}
    \label{fig:problem_with_0_eigenvalue_polar}
\end{figure}

A Lyapunov function $V(\vb{y})$ can be found to show that $(0,0,0)^T$ is a stable fixed point
\begin{equation}
\label{eq:liapunov}
    V(\vb{y}) = \frac{1}{2} (y_1^2 + y_2^2 + y_3^2) ,
\end{equation}
see appendix \ref{app:slow_mode_Problem} for details. Especially, we can show that the trajectories converge to the steady state.

\subsubsection{Amplitude-phase formulation for problem with $0$-eigenvalue}

We apply the theory from \cref{subsec:general_Weakly_nonlinear_system} and write the solution of the  problem with $0$-eigenvalue in the form of equation \eqref{eq:general_WKB_ansatz}:
\begin{align}
     \vb{y} (t)
     & = a_0(t) \begin{pmatrix} 0 \\ 0 \\1 \end{pmatrix}
 +  a_1(t) e^{\frac{i}{\epsilon}(t+ \epsilon \psi_1(t))} \frac{1}{\sqrt{2}} \begin{pmatrix} 1 \\ -i \\0 \end{pmatrix} + a_1(t) e^{-\frac{i}{\epsilon}(t+ \epsilon \psi_1(t))} \frac{1}{\sqrt{2}} \begin{pmatrix} 1 \\ i \\0 \end{pmatrix} 
 \\ 
 &=  \begin{pmatrix} \sqrt{2} a_1(t) \cos(\frac{t}{\epsilon} + \psi_1(t)) \\ \sqrt{2} a_1(t) \sin(\frac{t}{\epsilon} + \psi_1(t)) \\ a_0(t)\end{pmatrix} \label{eq:zeroEof_newCoords}
\end{align}


\subsubsection{Differential equations in the new coordinates}

We can show that the problem with $0$-eigenvalue admits the following form in the $a_0, a_1, \phi_1$-coordinates
\begin{align}
    \dot a_1 &= -\lambda (2 a_1^2 + a_0) a_1\\
    \dot \psi_1 &= 0 \\
    \dot a_0 &= -\lambda (2 a_1^2 + a_0) a_0 ,
    \label{eq:0eigenvalue_slow_system}
\end{align}
see appendix \ref{app:0eig_new_coordinates}.
The transformed equations contain no explicit $\frac{t}{\epsilon}$ terms on the right-hand side; in particular, the small parameter $\epsilon$ is no longer present there. This suggests that the transformed system evolves only on a slow time scale and, in particular, does not exhibit the fast oscillations observed in the original variables $y_1$, $y_2$, and $y_3$. \Cref{fig:problem_with_0_eigenvalue_polar} also indicates that the solution evolves on a slow time scale in the new coordinates.


Assume that $a_0(t) > 0$ and $a_1(t) > 0$. Then, we can show that the amplitude functions differ by a constant multiplicative factor. 
More precisely, we have
\begin{equation}
    \frac{a_1(t)}{a_0(t)} = \frac{a_1(t_0)}{a_0(t_0)} \equiv \mathrm{const.},
\end{equation}
 see appendix \ref{app:relation_a_a_0} for details.

\subsection{System 2: van der Pol oscillator}
\label{subsec:vdp_analysis}

The van der Pol oscillator is a second order ordinary differential equation given by
\begin{equation}
    \label{eq:vdP}
    \dv[2]{y}{t} + y = \epsilon (1-y^2) \dv{y}{t} ,
\end{equation}
It can be rewritten as a first order system 
\begin{equation}
\label{eq:vdP_system}
\begin{split}
    \dot y_1 &= y_2 \\
    \dot y_2 &= -y_1 + \epsilon (1-y_1^2)  y_2 .
    \end{split}
\end{equation}

We consider the case $0 < \epsilon \ll 1 $. Especially, the system has an attracting  limit cycle, and an unstable fixed point at $\vb{y}^{\star}=(0,0)^T$, \textit{i.\ e.} the right-hand side of \eqref{eq:vdP_system} vanishes at $\vb{y}^{\star}$. 
After applying the scaling of the time coordinate $\tau = \epsilon t$, the van der Pol oscillator admits the form of the weakly nonlinear system \eqref{eq:general_problem}. 
\footnote{In particular, we introduce $y^{\star}(\tau) = y^{\star}(\epsilon t) = y(t)$, reformulate \eqref{eq:vdP_system} in terms of $y^{\star}$, and finally drop the asterisks to avoid cumbersome notation.}
The matrix $\mathcal{L}$ and the nonlinear operator $\mathcal{N} (\vb{y})$ are given by
\begin{equation}
\label{eq:general_problem_vdp}
    \mathcal{L} = \begin{pmatrix}
        0 & 1 \\ -1 & 0
    \end{pmatrix} , \qquad
    \mathcal{N} (\vb{y}) = \begin{pmatrix}
        0 \\ g(\vb{y})  
    \end{pmatrix} , \qquad g(\vb{y}) =  (1- y_1^2) y_2.
\end{equation}
In the remainder of the paper, we will consider the rescaled equations.

The eigenvalues of the linear operator $\mathcal{L}$ are given by $\mu_{1} = i$ and $\overline{\mu_{1}} = -i$. The corresponding normalized eigenvectors are
$\vb{r}_{1} = \frac{1}{\sqrt{2} }
\begin{pmatrix} 1 \\ i \end{pmatrix}$ and $ \overline{\vb{r}}_{1} = \frac{1}{\sqrt{2} }
\begin{pmatrix} 1 \\ -i \end{pmatrix}$, and will be used below.  

The van der Pol system has two time scales, one due to the fast linear operator and slow one due to the slow nonlinear term. 

\subsubsection{Amplitude-phase formulation for the van der Pol oscillator}
\label{subsec:WKB_van_der_Pol}

The theory developed in \cref{subsec:general_Weakly_nonlinear_system} is applied to the rescaled van der Pol oscillator \eqref{eq:general_problem_vdp}. As the system is two dimensional and we have only one pair of complex-conjugate eigenvalues, we drop the index used in equation \eqref{eq:genaral_lin_sol}. In agreement with \eqref{eq:genaral_lin_sol}, the solution of the linear system can be written as
\begin{equation}
\label{eq:lin_sol_2}
    \vb{y}_{\text{lin}} = a e^{\frac{i}{\epsilon}(t+\epsilon \psi)} \vb{r} + a e^{-\frac{i}{\epsilon}(t+\epsilon \psi)} \overline{\vb{r}} .
\end{equation}

Moreover, using the results from \cref{subsec:general_Weakly_nonlinear_system}, we obtain the following representation of the solution of the weakly nonlinear system
 \begin{equation}
 \label{eq:WKB_ansatz_van_der_Pol}
 \begin{split}
     \vb{y} (t)
     &= a(t) e^{\frac{i}{\epsilon}\phi (t)} \vb{r} +  a(t) e^{-\frac{i}{\epsilon}\phi (t)} \overline{\vb{r}}  \\
     &= a(t) e^{\frac{i}{\epsilon}\phi (t)} \frac{1}{\sqrt{2}} \begin{pmatrix}  1 \\i  \end{pmatrix} +  a(t) e^{-\frac{i}{\epsilon}\phi (t)} \frac{1}{\sqrt{2}} \begin{pmatrix}
         1 \\ -i
     \end{pmatrix}\\
     &= \sqrt{2} a(t) \begin{pmatrix}
         \cos(\frac{\phi}{\epsilon}) \\
          - \sin(\frac{\phi}{\epsilon})  
     \end{pmatrix} ,
     \end{split}
 \end{equation}
 see equation \eqref{eq:general_WKB_ansatz} for comparison. Note that $\phi$ satisfies equation \eqref{eq:def_phi}.

\subsubsection{Equations for amplitude and phase}

The functions describing the amplitude $a(t)$ and the phase  $\phi(t)$ in  \eqref{eq:WKB_ansatz_van_der_Pol} obey the  relation (see appendix \ref{subsec:Amplitude_Phase_Equs} and see appendix \ref{subsec:furhter_details_rhs})
\begin{align}
\label{eq:wkb_amplitude}
    \dot a(t) &= \frac{1}{2} \begin{pmatrix}
        0 \\ g(\vb{y}) 
    \end{pmatrix}  \cdot \left( e^{-\frac{i}{\epsilon}\phi} \overline{\vb{r}}  + e^{\frac{i}{\epsilon}\phi} \vb{r}  \right)\\
    \dot \phi(t) &= 1 - \frac{i \epsilon}{2 a}  \begin{pmatrix}
        0 
        \\ g(\vb{y}) 
    \end{pmatrix}  \cdot \left( e^{-\frac{i}{\epsilon}\phi} \overline{\vb{r}}  - e^{\frac{i}{\epsilon}\phi} \vb{r}  \right) . \label{eq:wkb_phase}
\end{align}

\begin{figure}[h]
    \centering
  \includegraphics[width=0.75\linewidth]{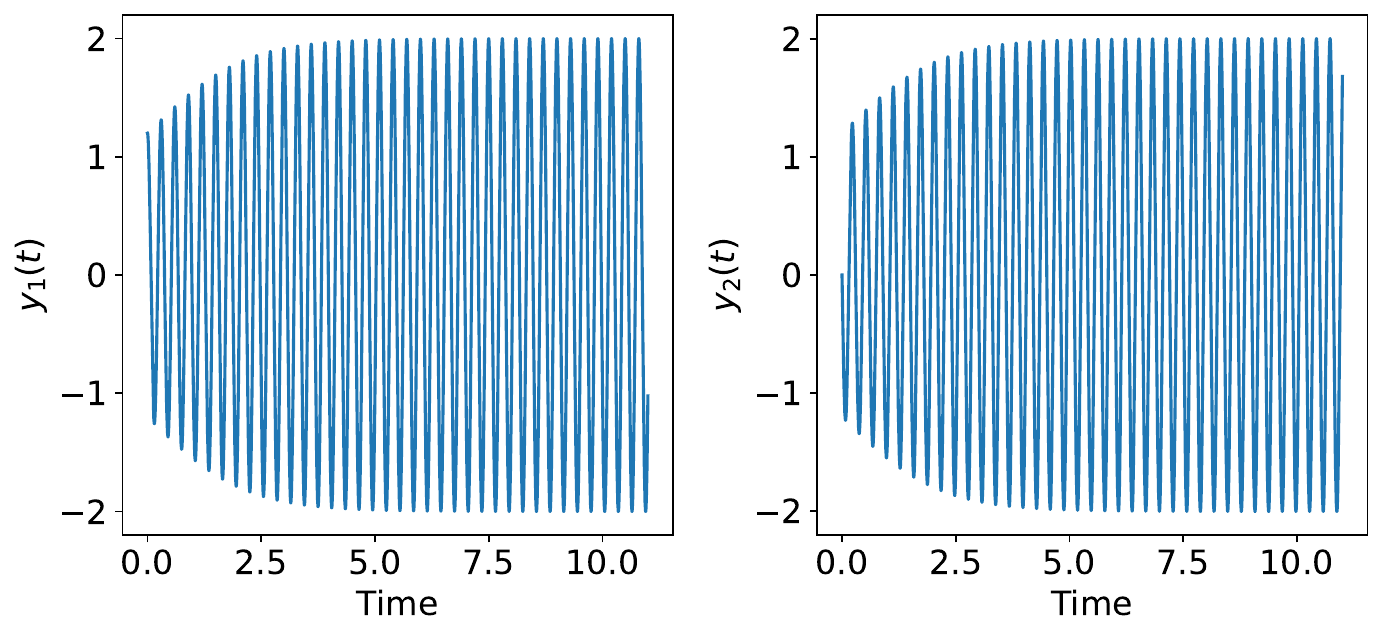}
    \caption{Numerical solution of the van der Pol system \eqref{eq:general_problem_vdp}. The original coordinates $y_1(t)$ and $y_2 (t)$ are shown for the case $\epsilon = \frac{0.3}{2 \pi}$.}
    \label{fig:vdp-y1_y2}
\end{figure}

\begin{figure}[h]
    \centering
    \includegraphics[width=0.75\linewidth]{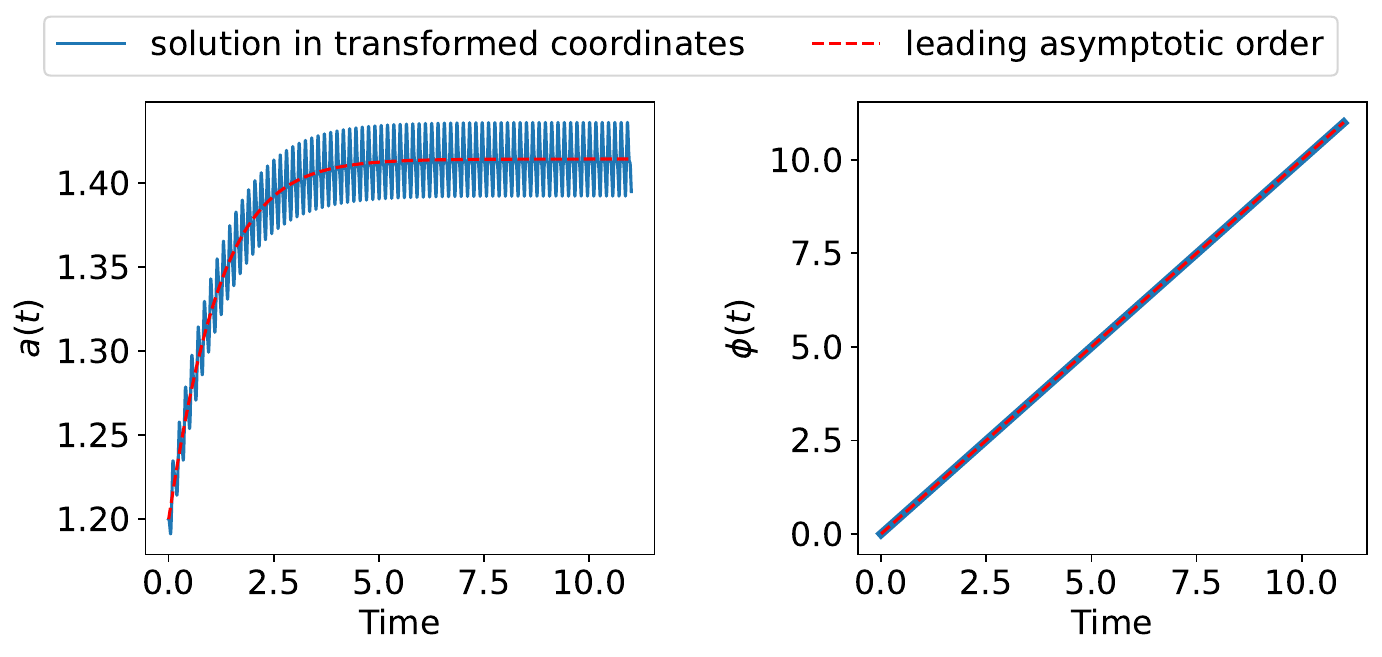}
    \caption{Numerical solution of the van der Pol system \eqref{eq:general_problem_vdp}. The polar variables $a(t)$ and $\phi(t)$ are shown for $\epsilon = \frac{0.3}{2\pi}$, where $\phi(t)$ is the shown phase variable and the actual angle is $\phi(t)/\epsilon$. For comparison, the leading order asymptotic solution of \eqref{eq:vdp_asymp_a} and \eqref{eq:vdp_asymp_psi}, in the form $t+ \epsilon \psi^{(0)}$, are shown .
}
    \label{fig:vdp-amplitude_phase}
\end{figure}

\cref{fig:vdp-y1_y2} and \cref{fig:vdp-amplitude_phase} show solutions to the  van der Pol system in the original coordinates $y_1, y_2$ and polar variables $a, \phi$. In particular, in the original coordinates, the rapid oscillations have amplitudes of order one, whereas in polar variables the leading-order dynamics are slow. Moreover, in the new coordinates, the rapid oscillations have smaller amplitudes and appear only at higher orders in $\epsilon$.
The appearance of rapid oscillations at higher orders is a crucial difference from the results for the problem with $0$-eigenvalue, where the fast oscillations are completely eliminated by the transformation \eqref{eq:general_WKB_ansatz}.

\subsubsection{Asymptotics of van der Pol system}

Motivated by the previous section, we assume that $a$ and $\psi$ depend on two time scales
\begin{equation}
    a(t) = a\left(t, \frac{t}{\epsilon} \right) = a(t,\tau) , \qquad \qquad \psi(t) = \psi\left(t, \frac{t}{\epsilon} \right) = \psi(t, \tau) .
\end{equation}
and employ the asymptotic expansions:
\begin{align}
    a(t) &= a^{(0)}(t) + \epsilon a^{(1)}(t, \tau) + \epsilon^2 a^{(2)}(t, \tau) + o(\epsilon^2) \\
    \psi(t) &= \psi^{(0)}(t) + \epsilon \psi^{(1)}(t, \tau) + \epsilon^2 \psi^{(2)}(t, \tau) + o(\epsilon^2) 
\end{align}

The leading order system is given by

\begin{align}
    a_t^{(0)} &= \frac{1}{2} a^{(0)} - \frac{1}{4} {a^{(0)}}^3, \label{eq:vdp_asymp_a}\\
    \psi_t^{(0)} &= 0 , \label{eq:vdp_asymp_psi}
\end{align}

for details see appendix \ref{app:Asymptotics_for_Van_der_Pol}.

\section{WKB-related time-stepping scheme}
\label{sec:WKB-related_time-stepping scheme}

In this section, we describe a novel time-stepping scheme that is adapted to problems of the type \eqref{eq:general_problem} and makes use of the transformation \eqref{eq:general_WKB_ansatz} without explicitly employing the governing equations in the transformed coordinates. The new method is tested on the examples presented in \cref{sec:Example_systems}. We now describe the building blocks of time-stepping scheme. In addition, \cref{alg:wkb-0-eigenvalue,alg:wkb-vdp} illustrate the structure of the WKB time-stepping scheme for the zero-eigenvalue problem and the van der Pol system. The detailed functions called within these algorithms are presented as pseudocode in appendix \ref{app:pseudocode}.

\subsection{Short time simulations}

First, short simulations of problem \eqref{eq:general_problem} are computed using a standard time integrator over intervals of length~$L$, \textit{i.\ e.}, $\mathcal{J}_n = [n \Delta T, n \Delta T + L]$, where $\Delta T$ is the size of the large time step and the index $n$ counts the the large time steps. Hence, the starting points of consecutive intervals are separated by $\Delta T$. Within an interval $\mathcal{J}_n$, the system undergoes only a few oscillations. The standard time integrator is applied with a small time step~$\Delta t$ that resolves the fast oscillations.
In this study, the classical fourth-order Runge–Kutta method (RK4) is used as a standard integrator. However, the approach is not limited to RK4, and other time integrators can be employed as well.

Ideally, the short time simulations are over a small number of oscillations. As the period of one oscillation is of order $O(\epsilon)$, the length of the short time window is of order $O(N\epsilon) = O(\epsilon)$, where $N$ denotes the number of oscillations over which we integrate. (We use the period of the oscillations in the linear system as a reference.)

\subsection{Determination of amplitude and phase from numerical data}
\label{subsec:Determination_of_Data}

The solution data from the short time simulations is given in original (Cartesian) coordinates. In order to work with the phase and amplitude formulation introduced through \eqref{eq:general_WKB_ansatz}, we extract the corresponding amplitude and phase variables from the data. The following step is a crucial ingredient in the construction of the time-stepping scheme, as it allows us to determine the amplitudes and phases by solving linear optimization problems and thereby avoids the need to solve  nonlinear optimization problems.

Suppose we compute a short time simulation over few oscillations in an interval of length $L$. This provides the data $y_{1,i}, y_{2,i}$ in the case of the van der Pol oscillator. The first index denotes the component of the two dimensional system \eqref{eq:general_problem}, \eqref{eq:general_problem_vdp}, whereas the second index $i$ refers to the grid point in the short time window of length $L$. To extract the information of the amplitude, the following relation is used
\begin{align}
\label{eq:data_a}
    a_i &= \frac{1}{\sqrt{2}} \sqrt{y_{1,i}^2 + y_{2,i}^2} 
\end{align}
see equation \eqref{eq:WKB_ansatz_van_der_Pol} for comparison. To find the data of the phase function, we exploit the equation
\begin{equation}
\label{eq:data_phi}
    \tan(\frac{\phi_i}{\epsilon})  = -  \frac{y_{2,i}}{y_{1,i}} .
\end{equation}
 In the case of the van der Pol oscillator, the data  $a_i$ and $\phi_i$ will be averaged, see \cref{subsec:Averaging_of_data}, and then used in the fitting method described in the \cref{subsec:Fitting_of_the_parameters}. 

 An analogous approach  applies to the problem with a $0$-eigenvalue, see appendix \ref{app:recover_0eig_ampl_phases}. However, since the system has no rapid oscillations in the transformed coordinates, smoothing operations can be neglected.
In general, the equations for computing the amplitude and phase functions, $a_j$ and $\phi_j$, must be adapted to the eigenvalue and eigenvector decomposition of the linear operator.

\subsection{Smoothing of the amplitude and phase data}
\label{subsec:Averaging_of_data}

To prevent the remaining fast oscillations from influencing the approximation of $a(t)$ and $\phi(t)$, where in the present study the approximations of $a(t)$ and $\phi(t)$ model the behavior on the slow time scale only, we define a smoothed approximation via a convolution-type integral. In this approximation, the remaining oscillations appearing at higher asymptotic orders are attenuated. We first describe the details of the smoothing operation applied in the present study. At the end, of this \cref{subsec:Averaging_of_data}, we illustrate how this procedure affects the approximations of the slow dynamics.

For the smoothing method, the following integral is evaluated
\begin{equation}
\label{eq:av_a_analytical}
    \overline{a}(t) = \frac{1}{\eta} \int_{- \frac{\eta}{2}}^{\frac{\eta}{2}} \rho \left( \frac{s}{\eta}  \right) a(t+s) \dd{s} ,
\end{equation}
where $\eta$ is denoted as the smoothing window. Typically, to attenuate the oscillations, $\eta$ spans the duration of a few periods. The function $\rho$ denotes the kernel function, defined by
\begin{equation}
\rho(s) =
\begin{cases}
C  \exp\Big(\frac{1}{(s-1/2)(s+1/2)}\Big), & -\frac{1}{2} < s < \frac{1}{2}, \\
0, & \text{otherwise.}
\end{cases}
\end{equation}
The constant $C$ is chosen such that $\|\rho\|_{L^1(\mathbb{R})} = 1$.
An analogous integral is formulated for the function in the phase $\phi$. Its averaged version is denoted as $\overline{\phi}$.

For the numerical evaluation of the integrals, we apply numerical quadrature and use the data $a_i$ and $\phi_i$, see \cref{subsec:Determination_of_Data}. In this study, we apply the trapezium rule

\begin{equation}
    \overline{a}_j = \frac{1}{\eta} \sum_{i \in \mathcal{I}} 
    \rho \left( \frac{s_i}{\eta} \right) {a}_{j+i}
 \Delta t ,
\end{equation}
where $\{s_i | \ i \in \mathcal{I}\} = \{ -\frac{\eta}{2}, -\frac{\eta}{2}+ \Delta t, \dots,  \frac{\eta}{2} - \Delta t, \frac{\eta}{2}  \}$, that is we introduce a discrete grid with steps size $\Delta t$ in the interval $[  -\frac{\eta}{2},  \frac{\eta}{2} ]$. For simplicity, we take $\eta$ to be an integer multiple of $\Delta t$, but the approach can be straightforwardly extended to the general case. Since
$\rho(\pm \tfrac12)=0$, the endpoint contributions in the trapezium rule vanish, and therefore no additional factor $1/2$ appears in the above expression. We evaluate an analogous discretization to find the $\overline{\phi}_j$. 



\begin{figure}[H]
    \centering
    \includegraphics[width=0.5\linewidth]{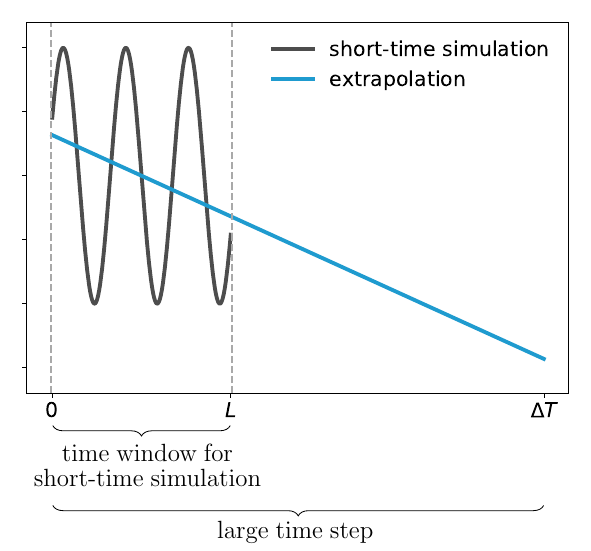}
    \caption{Sketch  of short time simulation in the interval $[0,L]$ and extrapolated fitted straight line: at the beginning, most data points of the short time simulation lie above the slow trend, while towards the end more points lie below it. This asymmetry biases the linear fit downward, leading to an overly negative slope.}
    \label{fig:sketch}
\end{figure}

\cref{fig:sketch} shows a sketch of a short time simulation of a highly oscillatory function over a small time interval. We want to approximate the slow (average) dynamics of this function by fitting a straight line to the data in this short time window.
However, this linear fit can be distorted by the fast oscillations: the fitted line is meant to serve as an approximation of the slow dynamics, but its slope may be biased by an uneven distribution of data points around the mean behavior. For example, if at the beginning of the time window most data points lie above the average trend and towards the end most points lie below it (as sketched in the figure), then the fitted line will be tilted downward more strongly than the true slow trend.
The opposite effect occurs if, at the beginning of the interval, most points lie below the average behavior and later above it: in that case, the fitted slope will be biased upward. In both situations, the fast oscillations distort the estimate of the slow dynamics obtained from the linear fit.

\begin{figure}[H]
    \centering
    \begin{minipage}{0.48\linewidth}
        \centering
        \includegraphics[width=\linewidth]{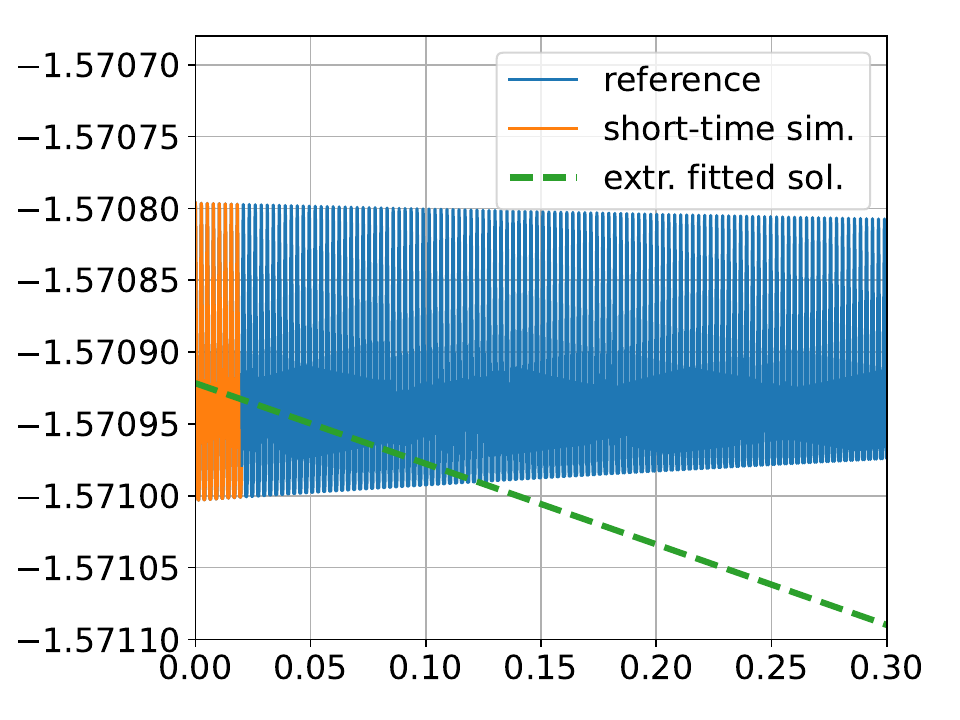}
        \caption{Approximation of $\psi (t)$ for the van der Pol system with initial value $\vb{y}_0 = (0.0, 1.0)^T$ and small parameter $\epsilon = \frac{5.2 \times 10^{-3}}{2 \pi}$. The blue line shows the reference solution, the orange line is the short time simulation in the first time window of length  $L= 0.02$ with short time step $\Delta t = 1.0 \times 10^{-6}$. The green line shows the extrapolated approximation after applying the fitting method to the data from the short time simulation without smoothing.} 
        \label{fig:psi_no_averaging}
    \end{minipage}
    \hfill
    \begin{minipage}{0.48\linewidth}
        \centering
        \includegraphics[width=\linewidth]{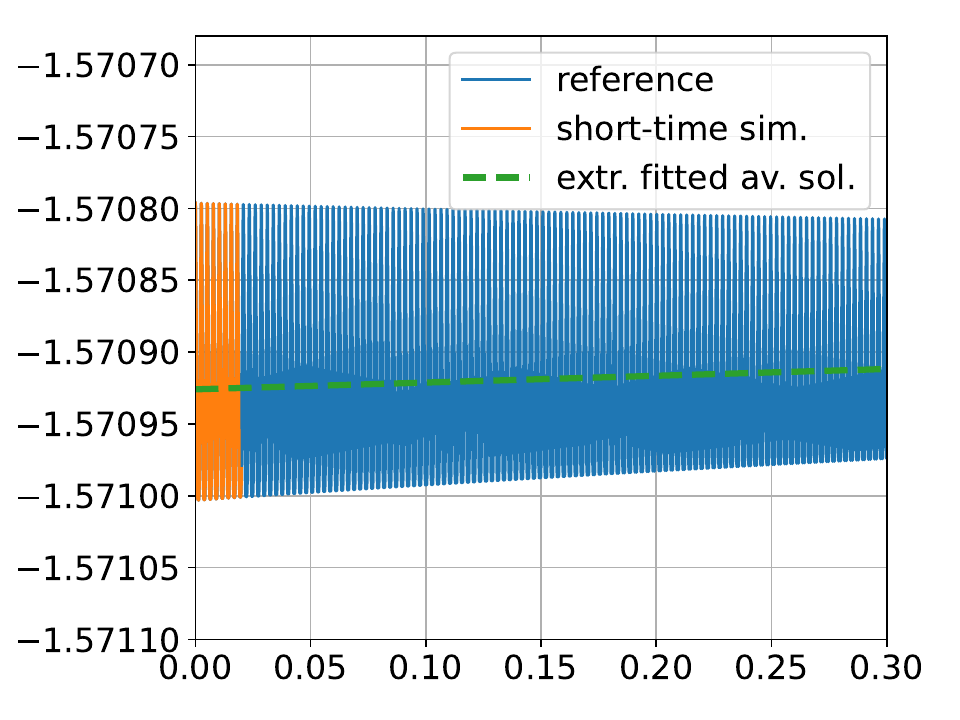}
        \caption{Approximation of $\psi (t)$ for the van der Pol system with initial value $\vb{y}_0 = (0.0, 1.0)^T$ and small parameter $\epsilon = \frac{5.2 \times 10^{-3}}{2 \pi}$. The blue line shows the reference solution, the orange line is the short time simulation in the first time window of length  $L= 0.02$ with short time step $\Delta t = 1.0 \times 10^{-6}$. The green line shows the extrapolated approximation after applying the fitting method to the data from the short time simulation with smoothing using a smoothing window $\eta = 1.0 \times 10^{-2}$.}
        \label{fig:psi_with_averaging}
    \end{minipage}
\end{figure}

This effect can be mitigated by applying a smoothing procedure to the data. In our concrete example of the van der Pol oscillator, we illustrate this for the coordinate $\psi$ in \cref{fig:psi_no_averaging,fig:psi_with_averaging}.
In \cref{fig:psi_no_averaging}, one observes exactly the distortion effect illustrated in the sketch: the local linear fit is biased by the fast oscillations, so that the fitted line is noticeably tilted and does not correctly represent the slow average behavior.
In contrast, \cref{fig:psi_with_averaging} shows how this effect is reduced when the smoothing procedure is used. The fit is much less biased by the oscillations but instead yields a straight line that stays roughly above the oscillations and tracks the slow trend. This remains true even when the fit is extrapolated by a coarse time step  to the next short time interval.

The study \cite{Wingate_Rosemeier_Haut_2023} investigates how smoothing operations mitigate fast oscillations for appropriately chosen smoothing windows. In addition, \cite{HautWingate2014} provides an estimate quantifying the mitigation of rapid oscillations (see equation (3.9) therein).


\subsection{Fitting of the parameters}
\label{subsec:Fitting_of_the_parameters}

Since in the general case, exact analytical descriptions for $a(t)$ and $\phi(t)$ 
are not available, we propose a method to compute numerical approximations. 
We will study the following ansatz for the WKB-based time-stepping methods applied to the van der Pol oscillator:
\begin{align}
\label{eq:polynomial_linear_ansatz}
    a_{\text{fit}}(t) &= a_{0,\mathrm{f}} + t a_{1,\mathrm{f}}  , &
    \phi_{\text{fit}}(t) &= \phi_{0,\mathrm{f}} + t \phi_{1,\mathrm{f}} . 
\end{align}
The parameters $a_{0,\mathrm{f}}, a_{1,\mathrm{f}}, \phi_{0,\mathrm{f}}, \phi_{1,\mathrm{f}}$ are chosen in such a way that the polynomial expansions   approximate the data $\overline{a}_j$ and $\overline{\phi}_j$ as accurately as possible, in particular we apply a least squares approximation. Replacing $a(t)$ and $\phi(t)$ by ansatz \eqref{eq:polynomial_linear_ansatz} in equation \eqref{eq:WKB_ansatz_van_der_Pol} for the van der Pol example, we obtain approximate solutions, introducing an error. This approximation is a superposition of wave modes with a slowly varying amplitude function $a_{\text{fit}}(t)$ and and a slowly varying function $\phi_{\text{fit}}(t)$, and is therefore related to a multi-modal WKB approximation.

For the problem with 0-eigenvalue, we employ an analogous ansatz, \textit{i. \ e.}\ a linear fit of $a_0, a_1, \phi_1$, and the same reasoning can be applied in this context.

\subsection{Extrapolation}

The WKB approximation is used to perform a large time step, $\Delta T$, and compute a numerical solution far beyond the short time window.  
First, we replace $a(t)$ and $\phi(t)$ in \eqref{eq:WKB_ansatz_van_der_Pol} in the case of the van der Pol oscillator by $a_{\text{fit}}$ and   $\phi_{\text{fit}}$. Then, we evaluate this relation for $t = \Delta T$.
This solution then serves as the starting point for the next time window; see \cref{fig:idea_WKB_method}, and we can perform the next short time simulation. Repeating this procedure allows us to advance the solution over  long time intervals.
 When we think in terms of a two-scale system, the small time-step, $\Delta t$, is used to resolve the fast dynamics, whereas the large time-step, $\Delta T$, is used to resolve the slow dynamics.

\begin{remark}[Connection to Euler-type discretization]
\label{rem:Connection_Euler-type_discretization}
To explain the first-order convergence observed in the following section, we provide a heuristic interpretation of the proposed scheme. On short time intervals of length $\mathcal{O}(\epsilon)$, the variables $a(t)$ and $\phi(t)$ describing the van der Pol oscillator are well-approximated by their first-order Taylor expansions. The linear least-squares fits used in the method can therefore be viewed as data-driven approximations of these local Taylor expansions. The subsequent extrapolation step is then equivalent to an explicit Euler update for the slow dynamics. This provides a heuristic explanation for the observed linear convergence rate with respect to the coarse time step $\Delta T$. Analogous arguments apply to the problem with $0$-eigenvalue.
\end{remark}

\begin{algorithm}[H]
\caption{WKB time-stepping for problem with 0-eigenvalue}
\label{alg:wkb-0-eigenvalue}
\begin{algorithmic}[1]
\Function{WKBSimulation1}{$num\_steps, L, \Delta T, \Delta t, x_0, \varepsilon, rhs = RHS, T_{min}$}
    \State $t_0 \gets Tmin$  \Comment{Initial time}
    \State $t_1 \gets t_0 + L$  \Comment{End of first short time window}
    
    \For{$n = 0$ \textbf{to} $num\_steps - 1$}
        \State \(\triangleright\) Short time simulation
        \State $short\_grid, [y_1, y_2, y_3] \gets$ 
        \Call{Integrate}{$rhs, x_0, t_0, t_1, \Delta t$}   

        \State \(\triangleright\) Amplitudes and phases
        \State $a_1\_short \gets$ \Call{FindAmplitude}{$y_1, y_2$}
        \State $\phi_1\_short \gets$ \Call{FindPhi}{$y_1, y_2, \varepsilon$}
        \State $a_0\_short \gets y_3$ 
        
        \If{$n > 0$}
            \State \(\triangleright\) Adjust phase with previous step
            \State $\phi\_short \gets \phi_1\_short - \phi_1\_short[0] + extrapolated\_\phi_1[-1]$  
        \EndIf

         \State \(\triangleright\)  Fit linear parameters
        \State $(slope\_a_1, intercept\_a_1) \gets$ linear\_fit$(short\_grid, a_1\_short)$
        \State $(slope\_\phi_1, intercept\_\phi_1) \gets$ linear\_fit$(short\_grid, \phi_1\_short)$
        \State $(slope\_a_0, intercept\_a_0) \gets$ linear\_fit$(short\_grid, a_0\_short)$ 
        
        \State $bridge\_grid \gets$ $[t_0, t_0 + \Delta T]$  \Comment{extrapolation grid}

        \State \(\triangleright\) Extrapolate amplitude and phase
        \State $extrapolated\_a_1 \gets slope\_a_1 \cdot bridge\_grid + intercept\_a_1$
        \State $extrapolated\_\phi_1 \gets slope\_\phi_1 \cdot  bridge\_grid + intercept\_\phi_1$
        \State $extrapolated\_a_0 \gets slope\_a_0 \cdot bridge\_grid + intercept\_a_0$ 

        \State \(\triangleright\) Inverse transform to original coordinates
        \State $y_1 \gets \sqrt{2} \cdot extrapolated\_a_1[-1] \cdot \cos\left(\frac{1}{\varepsilon} \cdot extrapolated\_\phi_1[-1]\right)$
        \State $y_2 \gets (-\sqrt{2}) \cdot extrapolated\_a_1[-1] \cdot \sin\left(\frac{1}{\varepsilon} \cdot extrapolated\_\phi_1[-1]\right)$
        \State $y_3 \gets extrapolated\_a_0[-1]$ 
        
        \State $x_0 \gets [y_1, y_2, y_3]$  \Comment{Update initial value for next step}
        \State $t_0 \gets bridge\_grid[-1]$
        \State $t_1 \gets t_0 + L$  \Comment{Update time window}
    \EndFor
    \State \Return $x_0$
\EndFunction
\end{algorithmic}
\end{algorithm}

\begin{algorithm}[H]
\caption{WKB time-stepping for van der Pol-system}
\label{alg:wkb-vdp}
\begin{algorithmic}[1]
\Function{WKBSimulation2}{$\varepsilon, x_0, \Delta t, L, num\_Steps, \Delta T, rhs=RHS\_vdP$}
 \State \(\triangleright\) Smoothing grid and kernel
    \State $grid\_, sc\_rho\_ \gets$ scaled\_rho$(4\pi \varepsilon, \Delta t)$
\State $\tilde L \gets len(grid\_)$

    \State \(\triangleright\) Loop over coarse time steps
    \State $new\_x_0 \gets x_0$
       \For{$n = 0$ \textbf{to} $num\_Steps-1$}
     \State \(\triangleright\) Short time simulation
        \State $short\_time\_grid, [y_1, y_2] \gets$
         \Call{Integrate}{$\text{rhs}, new\_x_0, n \Delta T, n \Delta T + L, \Delta t$}
        
        \State \(\triangleright\) Data of amplitude and phase
        \State $short\_a \gets$ \Call{FindAmplitude}{$(y_1,y_2)$}
        \State $short\_\phi \gets$ \Call{FindPhi}{$(y_1, -y_2)$}
        
        \If{$n > 0$}
        \State \(\triangleright\) Adjust phase with previous step
            \State $short\_\phi \gets short\_\phi - short\_\phi[0] + c + d \cdot n \cdot big\_dt$
        \EndIf

        \State \(\triangleright\) Smoothing of data
        \State $M \gets len(short\_time\_grid) - len(grid\_)$ \Comment{Number of averaged data points}
        \State Initialize $av\_grid, av\_\phi, av\_a$ as zero arrays of length $M$
        \State $idx0 \gets$ index of $\min(|grid\_|)$
        
        \For{$m = 0$ \textbf{to} $M-1$}
        \State \(\triangleright\) Trapezoidal rule for quadrature
            \State $av\_\phi[m] \gets \text{trapezoid}(sc\_rho\_ \cdot short\_\phi[m:\tilde L+m], x=grid\_)$
            \State $av\_a[m] \gets \text{trapezoid}(sc\_rho\_ \cdot short\_a[m:\tilde L+m], x=grid\_)$
            \State $av\_grid[m] \gets short\_time\_grid[idx0+m]$
        \EndFor

        \State \(\triangleright\) Parameter fitting
        \State $a, b \gets$ linear\_fit($av\_grid, av\_a$)
        \State $c, d \gets$ linear\_fit($av\_grid, av\_\phi$)
        
        \State \(\triangleright\) Extrapolation
        \State $bridge\_grid \gets$ $[n \cdot \Delta T, (n+1)\cdot \Delta T] $
        
        \State $y_1\_ext, y_2\_ext \gets$ \Call{Extrapolate}{$a, b, c, d, bridge\_grid$}
        \State $new\_iv \gets [y_1\_ext[-1], y_2\_ext[-1]]$
    \EndFor
\EndFunction
\end{algorithmic}
\end{algorithm}

\section{Numerical experiments}
\label{sec:Numerical_experiments}

\subsection{Numerical solution of the problem with $0$-eigenvalue}

In this section, we solve problem \eqref{eq:ODE_with_zero_frequency}, \eqref{eq:ODE_with_zero_frequency_2}. We perform tests for six values of $\epsilon$ spanning three orders of magnitude. 
For the coarse regime, we choose $\epsilon = \frac{1.0 \times 10^{-2}}{2 \pi}$ and $\epsilon = \frac{1.1 \times 10^{-2}}{2 \pi}$. The intermediate regime is given by $\epsilon = \frac{1.0 \times 10^{-3}}{2 \pi}$ and $\epsilon = \frac{1.1 \times 10^{-3}}{2 \pi}$. For the fine regime, we select $\epsilon = \frac{1.0 \times 10^{-4}}{2 \pi}$ and $\epsilon = \frac{1.1 \times 10^{-4}}{2 \pi}$.  Moreover, we use $\lambda = 1.0$. The initial value is given by $\vb{y}_0 = (\sqrt{2}/2, 0.0, \sqrt{2}/2)^T$. The problem is solved on the time interval $[0,5]$, \textit{i.\ e.} $T_{\mathrm{max}} =5$, that is if $N$ is the number of coarse time steps, the coarse time step size is $\Delta T = \frac{5}{N}$.

We solve the problem using the WKB-based scheme without smoothing the data, as described in Algorithm \ref{alg:wkb-0-eigenvalue}, since the dependence on $\epsilon$ is eliminated in the transformed coordinates, thereby removing the rapid oscillations. We compute a reference solution using the RK4 scheme. The time step sizes for the coarse, intermediate and fine regimes are $\Delta t_{\mathrm{ref}} = 1.0 \times 10^{-6}, 1.0 \times 10^{-7}, \text{ and } 1.0 \times 10^{-8}$, respectively.

The error is measured as a discrete mean $l_1$ error in a time window $[T_{\mathrm{max}}, T_{\mathrm{max}}+L]$ in this study. (The error at the final point does not seem to be representative for the slow dynamics.) 
Let $r_n$, $n=0,\dots,N$, denote a uniform grid on $[T_{\max}, T_{\max}+L]$.
We denote by $x_n$ the numerical solution of a quantity of interest and by $x_{n}^{\mathrm{ref}}$ the reference solution. The error at grid point $r_n$ is defined as
\begin{equation}
    e_n = x_n - x_{n}^{\mathrm{ref}} .
\end{equation}
We define the error in the interval $[T_{\max}, T_{\max}+L]$ as
\begin{equation}
    \| \vb{e} \|_{l_1} = \frac{1}{N+1} \sum_{n=0}^N |e_n| .
\end{equation}
The figures in this section show the error of the first component of the problem with $0$-eigenvalue.


In the tests we use $\Delta T \ge L$. This is a reasonable choice since in the other case the method would essentially reduce to the linear fit of the solution in the short time window described in \cref{subsec:Fitting_of_the_parameters}, thereby losing the benefits of the extrapolation.  

\begin{figure}[h]
    \centering

    \begin{minipage}{0.48\linewidth}
        \centering
        \includegraphics[width=\linewidth]{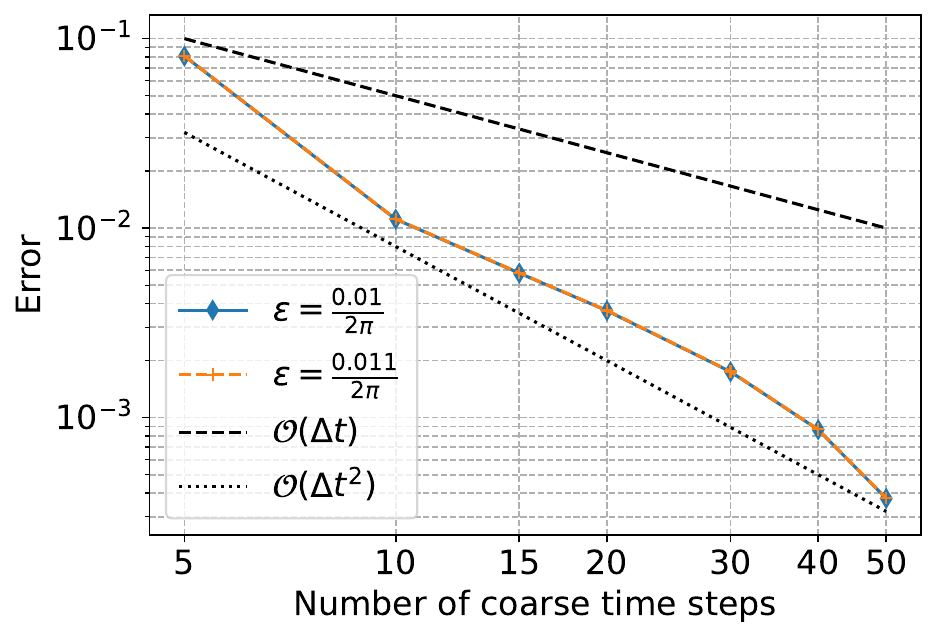}
        \caption{Numerical error for the problem with 0-eigenvalue plotted against the number of coarse time steps $N$ for two large values of $\epsilon$. The reference lines indicate first- and second-order convergence in $\Delta T$.}
        \label{fig:zero_eigen_large_eps}
    \end{minipage}
    \hfill
    \begin{minipage}{0.48\linewidth}
        \centering
        \includegraphics[width=\linewidth]{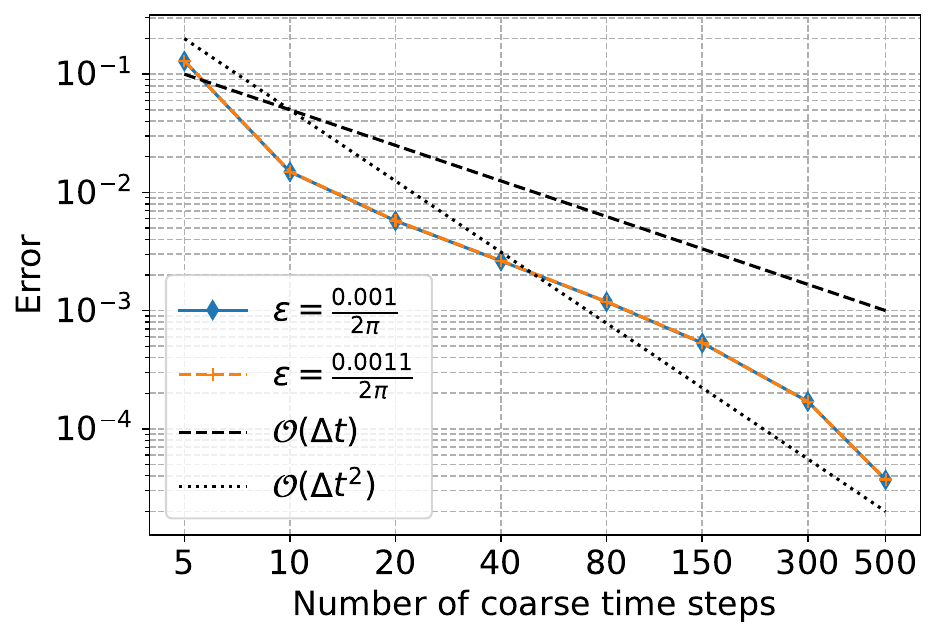}
        \caption{Numerical error for the problem with 0-eigenvalue plotted against the number of coarse time steps $N$ for two medium values of $\epsilon$. The reference lines indicate first- and second-order convergence in $\Delta T$.}
        \label{fig:zero_eigen_medium_eps}
    \end{minipage}

\end{figure}

\cref{fig:zero_eigen_large_eps} shows the numerical errors for the larger values of $\epsilon$. For the WKB scheme, we use the short time step $\Delta t = 1.0 \times 10 ^{-5}$ and the length of the short time window $L = 0.1$. To guarantee that the coarse time step is not shorter than the length of the short time windows $L$, we do not use more than 50 coarse time steps.  
We observe that the error decreases as $\Delta T$ is reduced. The experimental convergence order is roughly between 1 and 2. The difference between the two values of $\epsilon$ is not visually distinguishable. We assume that the error does not originate from the short time simulation, since $\Delta t$ is sufficiently small. This suggests that the dominant source of error lies in the fitting and extrapolation procedure in the transformed coordinates, where the underlying equation is independent of $\epsilon$. In contrast, we will see that for the van der Pol system, where the transformed equation does depend on 
$\epsilon$, the error curves do not coincide (at least not for small $\Delta T$).

In \cref{fig:zero_eigen_medium_eps}, we illustrate the errors for the medium values of $\epsilon$. The short time step is chosen as $\Delta t = 1.0 \times 10 ^{-6}$ and for the length of the short time window we use $L = 0.01$. To ensure that the coarse time step is not shorter than the length of the short time window $L$, we restrict the number of coarse time steps to at most 500.  We make the same observation as in the previous figure: the error is reduced when $\Delta T$ is decreased and the difference between the two values of $\epsilon$ is not visually distinguishable. We further observe approximately first-order convergence as $\Delta T$ is reduced.

\begin{figure}[H]
    \centering
    \includegraphics[width=0.6\linewidth]{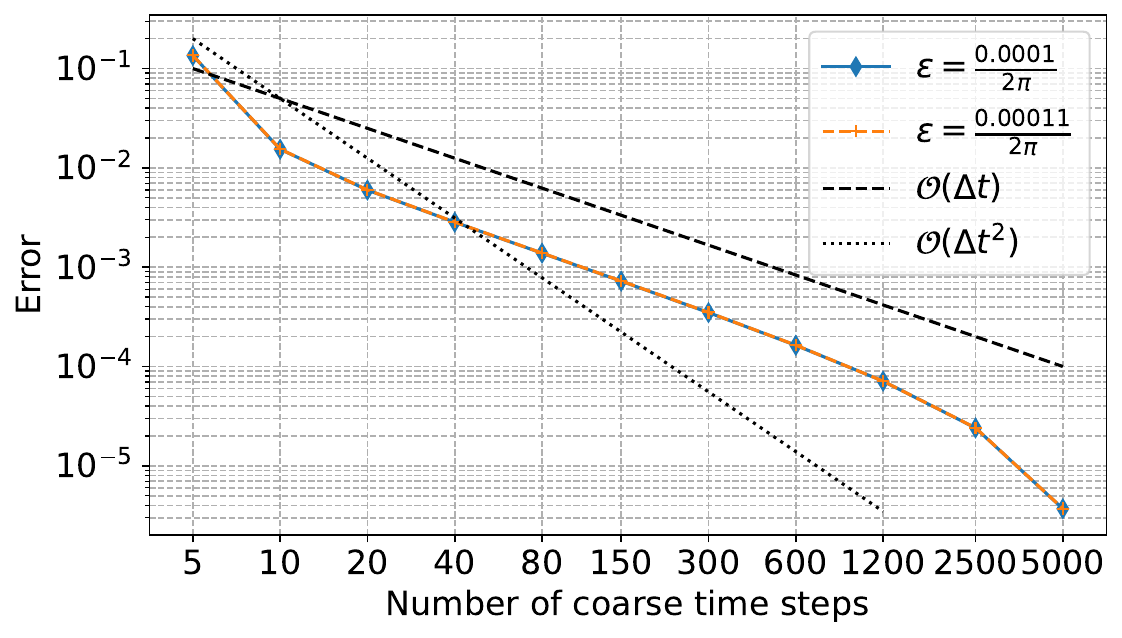}
    \caption{Numerical error for the problem with 0-eigenvalue plotted against the number of coarse time steps $N$ for two small values of $\epsilon$. The reference lines indicate first- and second-order convergence in $\Delta T$.}
    \label{fig:zero_eigen_small_eps}
\end{figure}

\Cref{fig:zero_eigen_small_eps} depicts the numerical error for the small values of $\epsilon$. We use the short time step $\Delta t = 1.0 \times 10 ^{-7}$ and the length of the short time window is $L = 0.001$ in the experiment. To guarantee that the coarse time step is not shorter than the length of the short time windows $L$, we restrict the number of coarse time steps to at most 5000. As in the previous figures, we observe that the error is reduced when $\Delta T$ is decreased and the difference between the two values of $\epsilon$ is not visually distinguishable. We further observe approximately first-order convergence as $\Delta T$ is reduced.

The figures indicate that we can achieve higher accuracies for smaller values of $\epsilon$, since smaller coarse time steps $\Delta T$ become possible. This is due to the fact that the short time window length $L$ decreases with $\epsilon$, allowing us to reduce $\Delta T$ without violating the constraint $\Delta T \ge L$, see \cref{fig:zero_eigen_large_eps,fig:zero_eigen_medium_eps,fig:zero_eigen_small_eps}.

In the figures, we observe at least first-order (linear) convergence, as expected. This can be explained by analogy to the explicit Euler method; see \textit{Remark} \ref{rem:Connection_Euler-type_discretization}.

Finally, we comment on the computational cost. Although the stiffness of the problem increases as $\epsilon$ decreases, the overall number of computational operations remains independent of $\epsilon$ under a fixed coarse discretization. In particular, if we fix the number of coarse time steps, \textit{e.g.} $N = 40$, the computational cost is independent of the choice of $\epsilon$. This is due to the fact that both the short time window length $L$ and the fine time step $\Delta t$ scale with $\epsilon$, so that the cost of the short time simulations adjusts accordingly.

\subsection{Numerical solution of the Van der Pol System}




\begin{figure}[h]
    \centering
    \begin{minipage}{0.48\linewidth}
        \centering
        \includegraphics[width=\linewidth]{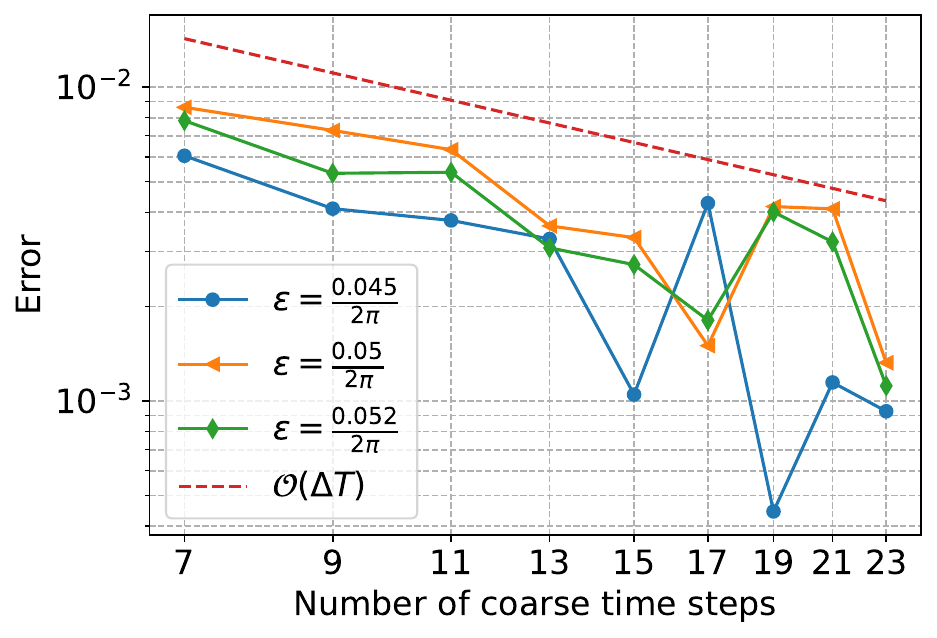}
        \caption{Error of the numerical approximation of the van der Pol oscillator for larger values of $\epsilon$. The dashed line serves as a reference for linear convergence order.}
        \label{fig:vdp_larger_eps_error}
    \end{minipage}
    \hfill
    \begin{minipage}{0.48\linewidth}
        \centering
        \includegraphics[width=\linewidth]{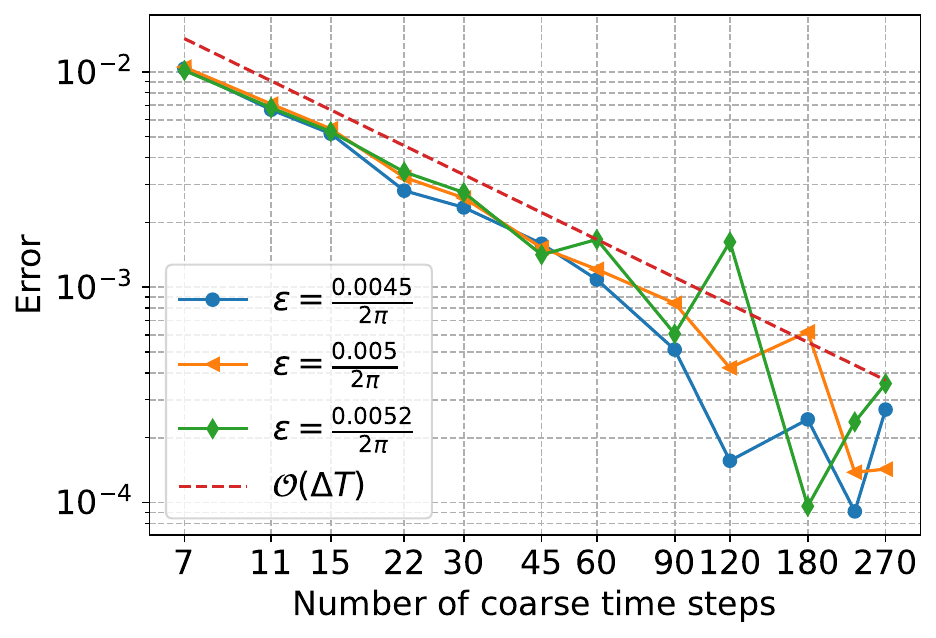}
        \caption{Error of the numerical approximation of the van der Pol oscillator for medium values of $\epsilon$. The dashed line serves as a reference for linear convergence order.}
        \label{fig:vdp_medium_eps_error}
    \end{minipage}
\end{figure}



\begin{figure}[H]
    \centering
    \includegraphics[width=0.6\linewidth]{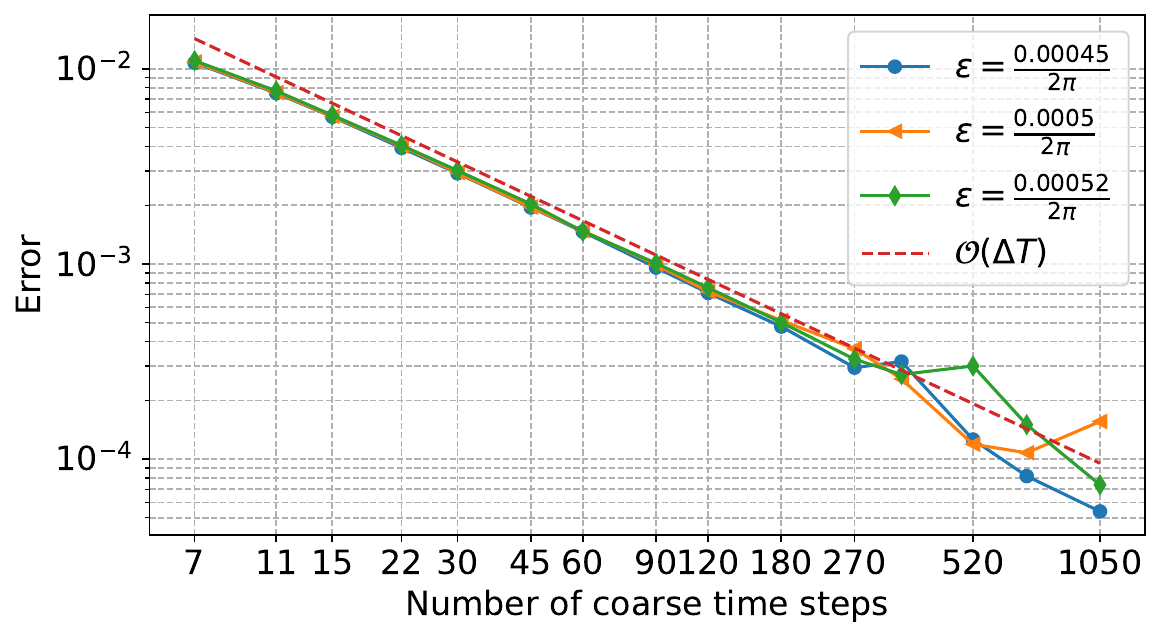}
    \caption{Error of the numerical approximation of the van der Pol oscillator for three small values of $\epsilon$. The dashed line serves as a reference for linear convergence order.}
    \label{fig:vdp_small_eps_error}
\end{figure}

We consider the van der Pol oscillator \eqref{eq:general_problem_vdp} on the time domain $[0,5]$, \textit{i. e.} $T_{\mathrm{max}} = 5$, a multiscale system also in the transformed coordinates, where the scale separation depends on $\epsilon$. 
The inital value is given by $\vb{y}_0 = (0.0, 1.0)^T$. 
Smaller values of $\epsilon$ correspond to stronger separation, making accurate approximation more challenging for standard numerical integrators.
In \cref{fig:vdp_larger_eps_error,fig:vdp_medium_eps_error,fig:vdp_small_eps_error}, we show the errors of the WKB-based scheme for three regimes of $\epsilon$ (large, medium, and small) as a function of the coarse time step $\Delta T$.
We computed reference solutions using the RK4 scheme with time step sizes $\Delta t_{\mathrm{ref}} = 1.0 \times 10^{-6}, 1.0 \times 10^{-7}, \text{ and } 1.0 \times 10^{-8}$ for the coarse, intermediate, and fine regimes, respectively.

In \cref{fig:vdp_larger_eps_error}, we show the error for three larger values of $\epsilon = \frac{4.5 \times 10^{-2}}{2\pi}, \frac{5.0 \times 10^{-2}}{2\pi}, \frac{5.2 \times 10^{-2}}{2\pi}$.  
For the short time simulations, a fine time step of $\Delta t = 1.0 \times 10^{-5}$ was chosen, and the length of the short time windows is $L = 2.0 \times 10^{-1}$. A smoothing window of $\eta = 1.0 \times 10^{-1}$ was used for temporal smoothing.
The difference between a reference solution with time step size $\Delta t_{\mathrm{ref}}$ and the solution with time step $\Delta t$ on intervals of length $L$ is below $10^{-8}$, indicating that the impact of the numerical error due to the short time simulations can be safely neglected. This error estimate remains valid for all regimes considered for the van der Pol system.

\Cref{fig:vdp_medium_eps_error} presents the results for three medium values of $\epsilon = \frac{4.5 \times 10^{-3}}{2\pi}, \frac{5.0 \times 10^{-3}}{2\pi}, \frac{5.2 \times 10^{-3}}{2\pi}$.  
The short time simulations were performed with a fine time step of $\Delta t = 1.0 \times 10^{-6}$, over  intervals of length $L = 2.0 \times 10^{-2}$. We used a smoothing window of $\eta = 1.0 \times 10^{-2}$ for temporal smoothing. 

Finally, \cref{fig:vdp_small_eps_error} shows the error for three small values of $\epsilon = \frac{4.5 \times 10^{-4}}{2\pi}, \frac{5.0 \times 10^{-4}}{2\pi}, \frac{5.2 \times 10^{-4}}{2\pi}$.
For the short time simulations, a time step of $\Delta t = 1.0 \times 10^{-7}$ was employed, and the system was integrated over short time windows of length $L = 2.0 \times 10^{-3}$. Temporal smoothing operations were applied using a smoothing window of $\eta = 1.0 \times 10^{-3}$.

Importantly, the computational cost of the method remains again independent of the choice of $\epsilon$, provided that the number of coarse time steps $\Delta T$ is kept constant. This is because both the fine time step and the length of the short time interval scale with $\epsilon$ in our approach across the regimes considered. As a result, for instance in the figures, each point corresponding to $\Delta T = 15$  coarse steps represents roughly the same computational effort, regardless of the value of $\epsilon$.

We approximate the slow evolution of the van der Pol oscillator's amplitude function $a(t)$ and phase function $\phi(t)$. 
These slow dynamics, however, are superimposed by the fast oscillations, whose amplitude is of order $\mathcal{O}(\epsilon)$ according to our investigations in \cref{subsec:vdp_analysis}. 
Consequently, one cannot generally expect to achieve an accuracy better than $\epsilon$, unless the approximations at the final time coincide by chance with the correct phase of the fast oscillation. 
Our fitting procedure is designed to capture the slow evolution of $a(t)$ and $\phi(t)$. Therefore, we expect that the fast oscillations limit the achievable accuracy to $\mathcal{O}(\epsilon)$.

In \cref{fig:vdp_larger_eps_error}, $\epsilon$ is roughly of order $10^{-2}$ (the absolute values are slightly smaller, but still of the same order). 
After only 7 coarse time steps, the accuracy of our method already exceeds $10^{-2}$, which explains why the error curves appear irregular: they are dominated by the neglected fast oscillations. 
For comparison, the dashed line shows the expected linear convergence versus the coarse time step $\Delta T$. Also, the fact that the error curves for different values of $\epsilon$ do not overlap indicates the impact of the fast oscillations on the error.

In \Cref{fig:vdp_medium_eps_error}, the error exhibits two distinct phases. 
In the first phase, the error decreases roughly linearly with the coarse time step $\Delta T$. 
In this phase, the error can be interpreted as being dominated by the approximation error of the slow evolution of $a(t)$ and $\phi(t)$.
It decreases roughly linearly with the coarse time step $\Delta T$, and can thus be reduced by taking smaller coarse steps.
In the first phase of the error curves, the results for different $\epsilon$ values almost coincide. 
This also indicates that the error is dominated by the approximation of the slow evolution of $a(t)$ and $\phi(t)$, rather than by the phase of the fast oscillations. 
After approximately 45 coarse steps, the error enters a second phase where it appears more irregular and reaches an order of magnitude of about $2 \times 10^{-3}$, which is consistent with the underlying $\epsilon$ values of roughly $10^{-3}$. The behavior in the second phase is similar to what is observed in \cref{fig:vdp_larger_eps_error}, where the error is dominated by the fast oscillations.

In \Cref{fig:vdp_small_eps_error}, a similar behavior to the medium-$\epsilon$ case is observed, with the $\epsilon$ values reduced by an additional factor of ten. 
Again, two distinct phases appear in the error curve. 
However, the first phase is significantly longer: the linear convergence persists up to approximately 270 coarse time steps, resulting in an error of about $3 \times 10^{-4}$. 
Beyond this point, the error starts to exhibit irregular behavior due to the fast oscillations.
This illustrates that smaller $\epsilon$ values allow for higher achievable accuracy in the approximation of the slow modulation.

Overall, the error is initially dominated by the slow dynamics, and only later influenced by fast oscillations, with the magnitude roughly scaling with $\epsilon$.





\section{Discussion and outlook}
\label{sec:Discussion_Outlook}

In this work, we have derived a transformation  for  differential equations describing oscillatory systems and used it to construct an associated time-stepping strategy based on amplitude–phase variables. The resulting method exhibits a close analogy to WKB-type approaches, as well as a structural similarity to classical explicit time integration schemes, especially the Euler method, when viewed on the level of the slow dynamics. Importantly, the evolution is carried out in the original (Cartesian) variables, and the transformed right-hand side does not need to be evaluated explicitly. 
The transformation itself encodes the multiscale structure of the system, leading to a formulation in which fast oscillations are either eliminated or strongly attenuated.
The method was tested for two systems of ordinary differential equations, including a discussion of the achieved accuracy.

A key property of the proposed approach is that its computational cost remains bounded as $\epsilon \to 0$, even in regimes where the underlying problem becomes increasingly stiff due to rapid oscillations. At the same time, the method exhibits intrinsic accuracy limitations: in particular, in the presence of higher-order terms, \textit{i.\ e.} terms of order $\mathcal{O}(\epsilon)$, in the transformed system that still contain fast oscillations, the error cannot, in general, be reduced below $\mathcal{O}(\epsilon)$ without further refinement of the transformation.

The present method bears some resemblance to the heterogeneous multiscale method (HMM) introduced in \cite{EEngquist2003} and applied to time-stepping problems, for example, in \cite{SharpTsaiEngquist2005}. 
In particular, short simulations on small time intervals are used to advance the solution on a larger time scale. 
In HMM, these simulations are typically employed to estimate effective macroscopic quantities, such as a coarse-grained right-hand side, often obtained via temporal averaging.
Instead, our approach relies on a transformation of the solution variables, and the micro-simulations are used to evolve the transformed system rather than to approximate a coarse-grained vector field.

At first sight, the present approach bears a structural resemblance to modulated Fourier expansions (MFE), a framework for the analysis of differential equations describing oscillatory systems \cite{HairerLubich2000,HairerLubichWanner2006}. 
In standard MFE, the solution is represented as a sum of Fourier modes with linear phases and slowly varying amplitudes.
The present method, in contrast, allows for more general phase functions, \textit{i.\ e.} we formulated the phase with a function $\psi(t)$, enabling nonlinearly varying oscillations. 
Moreover, our approach is based on the interpretation that the solution is represented in a set of polar coordinates, instead of a decomposition into oscillatory Fourier modes with slowly varying amplitudes.

In contrast to exponential integrators, see \cite{HochbruckLubich1997,CoxMatthews2002}, which rely on an exact treatment of the linear operator combined with a quadrature approximation of a nonlinear variation-of-constants integral, the present approach is based on a transformation of the solution variables and does not require an explicit representation of such an integral formulation.

In \cite{HautWingate2014,PeddleHautWingate2019,RosemeierHautWingate2024}, the authors develop Parareal methods for oscillatory problems. The Parareal method is a time-parallel algorithm that employs a coarse and a fine propagator. The present work should be contrasted with these coarse propagators: in the cited approaches, the governing equations are first transformed into a standard form, also denoted as modulation equation, \textit{i.\ e.}, a different transformation is applied. An HMM-type method is then applied to the transformed problem. In contrast, in the present method the untransformed, unaveraged equations are solved directly. Amplitude and phase coordinates are subsequently extracted from the solution data, and if smoothing is used, it is applied to these quantities. Finally, a fitting procedure approximates the polar variables as polynomials, which is then used for extrapolation over a large time step to the next time window.

The novel WKB-based method may be viewed as a candidate for a coarse propagator in Parareal-type algorithms for highly oscillatory problems, as it provides a good approximation of the oscillatory dynamics while the computational workload remains bounded even for small $\epsilon$. The full dynamics can still be resolved by the fine propagator, since the workload can be distributed across processors, making the approach potentially faster than a sequential fine-scale simulation.

An interesting direction for future research is the development of further transformations that reduce or eliminate the remaining $\epsilon$-dependent error induced by fast oscillations. In some systems, such as the problem with 0-eigenvalue, the present transformation is able to completely remove the $\epsilon$-dependence. This raises the question which structural properties of the underlying system enable such a complete elimination, and why in other cases, such as the van der Pol oscillator, fast oscillatory components persist after transformation.
Understanding these differences may provide insight into whether classes of oscillatory systems can be further transformed in such a way that the $\epsilon$-dependence is either entirely removed or shifted to higher-order terms.


A further promising avenue concerns the observed analogy between the proposed extrapolation procedure and the explicit Euler method. This suggests a deeper connection to classical time integration schemes and raises the question whether higher-order generalizations, such as Runge--Kutta or multistep-type methods, can be constructed within the same framework to improve accuracy on the slow time scale.

Another important extension is the application of the method to partial differential equations with highly oscillatory linear terms. Relevant examples include the rotating shallow water equations and the stratified Boussinesq equations. It would be interesting to investigate whether the present amplitude--phase based time-stepping strategy can be extended to yield efficient and robust schemes.

A further, more exploratory direction concerns extrapolation in $\epsilon$. The idea is whether a solution computed for a given value of $\epsilon$ can be used, via the present transformation framework, to construct approximations for nearby values of $\epsilon$ without recomputing the full time-stepping simulation. This would potentially allow for a more efficient reuse of existing simulation data, avoiding the need for repeated fine-scale computations.


\textbf{Acknowledgement:} This work was supported by the Collaborative Research Centre \emph{CRC 1114 Scaling Canscades in Complex Systems}, Project number 235221301, funded by Deutsche Forschungsgemeinschaft (DFG). The authors gratefully acknowledge the support. 
J.R. gratefully acknowledges support by the Deutsche Forschungsgemeinschaft
(DFG, German Research Foundation) under Germany's Excellence Strategy -- The
Berlin Mathematics Research Center MATH+ (EXC-2046/1, EXC-2046/2, project ID:
390685689).

\appendix

\section{Standard form of a weakly nonlinear problem}
\label{app:derivation_standard form}

We follow \cite{sanders_etal_2007} and consider
the perturbation problem 
\begin{equation}
\label{eq:perturbation_probelm}
    \dv{{\bf u}}{t} = A {\bf u} + \epsilon {\bf f} ({\bf u}, t, \epsilon),  \qquad \vb{u}(0) = \vb{u}_0
\end{equation}
where $A$ is a $n \times n$-matrix. The unperturbed problem takes the form
\begin{equation}
    \dv{{\bf v}}{t} = A {\bf v}, \qquad \vb{v}(0) = \vb{u}_0 .
\end{equation}

The solution of the linear problem is given by $\vb{v}(t) = e^{At} \vb{u}_0$.
To derive the standard form, variation of constants is applied, that is we replace $\vb{u}_0$ by a time-dependent function $\vb{w}(t)$.
\begin{equation}
    {\bf u}(t) = e^{At} {\bf w}(t)
\end{equation}
Thus, we have 
\begin{equation}
    \dv{\bf u}{t} = \dv{(e^{At})}{t} {\bf w}(t) + e^{At} (t) \dv{{\bf w}(t)}{t} = A  e^{At} {\bf w}(t) + \epsilon {\bf f} ( e^{At} {\bf w}(t), t, \epsilon) 
\end{equation}
We use that $e^{At}$ satisfies the differential equation $\dv{(e^{At})}{t}= A e^{At}$. Consequently, the standard form of \eqref{eq:perturbation_probelm} is given by
\begin{equation}
    \dv{\bf w}{t} = \epsilon e^{-At} {\bf f} ( e^{At} {\bf w}(t), t, \epsilon) . 
\end{equation}

\section{Details for the problem with $0$-eigenvalue}
\label{app:slow_mode_Problem}

\subsection{(Lyapunov-)stability of the steady state}

The origin $\vb{y}^{\star} = (0,0,0)^T$ is a steady state. Equation \eqref{eq:liapunov} satisfies the conditions of a Lyapunov function:
\begin{itemize}
    \item $V(\vb{y}^{\star}) = 0$
    \item \( V(\vb{y}) > 0\) for all other states \(\vb{y}\) in a neighborhood around $\vb{y}^{\star}$ 
    \item \(\dot{V}(\vb{y}) \le 0\), where 
     \(\dot{V}(\vb{y}) = y_1 \dot y_1 + y_2 \dot y_2 + y_3 \dot y_3 = -\lambda \| \vb{y} \|^4 \)
\end{itemize}

\subsection{Derivation of the ODE in the new coordinates}
\label{app:0eig_new_coordinates}

According to equation \eqref{eq:zeroEof_newCoords}, we have
\begin{equation}
    y_1 = \sqrt{2} a_1 (t) \cos(\frac{t}{\epsilon}+ \psi_1 (t)) \qquad
    y_2 = \sqrt{2} a_1 (t) \sin(\frac{t}{\epsilon}+ \psi_1 (t)) \qquad
    y_3 = a_0 (t)
\end{equation}

The temporal derivatives are expressed in the new coordinates:
\begin{align}
    \dot y_1 &= \sqrt{2} \dot a_1 \cos(\frac{t} {\epsilon}+\psi_1) - \sqrt{2} a_1 \left( \frac{1}{\epsilon} + \dot \psi_1 \right) \sin(\frac{t}{\epsilon}+\psi_1)\\
    \dot y_2 &= \sqrt{2} \dot a_1 \sin(\frac{t}{\epsilon}+\psi_1) + \sqrt{2} a_1 \left( \frac{1}{\epsilon} + \dot \psi_1 \right) \cos(\frac{t} {\epsilon}+\psi_1) \\
    \dot y_3 &= \dot a_0
\end{align}

The nonlinearity can also be rewritten in the new coordinates:

\begin{align}
    y_1^2 + y_2^2 + y_3^2 &= 2a_1^2 \left(\cos(\frac{t}{\epsilon}+\psi_1(t))\right)^2 + 2a_1^2 \left(\sin(\frac{t}{\epsilon}+\psi_1(t))\right)^2 + a_0^2 \\
    & = 2a_1^2 + a_0^2 
\end{align}

For the right-hand side, we find
\begin{align}
     (y_1^2 + y_2^2 + y_3^2) y_1 &= (2a_1^2 + a_0^2) \sqrt{2} a_1 \cos(\frac{t}{\epsilon}+\psi_1) \\
     (y_1^2 + y_2^2 + y_3^2) y_2 &= (2a_1^2 + a_0^2) \sqrt{2} a_1 \sin(\frac{t}{\epsilon}+\psi_1)  \\
     (y_1^2 + y_2^2 + y_3^2) y_3 &= (2a_1^2 + a_0^2) a_0
\end{align}

Thus, the differential equation rewritten in the new coordinates:

\begin{align}
    \sqrt{2} \dot a_1 \cos(\frac{t} {\epsilon}+\psi_1) - \sqrt{2} a_1 \left( \frac{1}{\epsilon} + \dot \psi_1 \right) \sin(\frac{t}{\epsilon}+\psi_1) &= \frac{\sqrt{2}}{\epsilon} a_1 \sin(\frac{t}{\epsilon}+\psi_1) - \lambda (2a_1^2 + a_0^2) \sqrt{2} a_1 \cos(\frac{t}{\epsilon}+\psi_1) \label{eq:zMode_newC_1}\\
    \sqrt{2} \dot a_1 \sin(\frac{t}{\epsilon}+\psi_1) + \sqrt{2} a_1 \left( \frac{1}{\epsilon} + \dot \psi_1 \right) \cos(\frac{t} {\epsilon}+\psi_1) &= -\frac{\sqrt{2}}{\epsilon} a_1 \cos(\frac{t}{\epsilon}+\psi_1) - \lambda (2a_1^2 + a_0^2) \sqrt{2} a_1 \sin(\frac{t}{\epsilon}+\psi_1)  \label{eq:zMode_newC_2}\\
     \dot a_0 &= - \lambda (2a_1^2 + a_0^2) a_0 \label{eq:zMode_newC_3}
\end{align}

Equation \eqref{eq:zMode_newC_3} describes the temporal evolution of $a_0$.
Computing $\eqref{eq:zMode_newC_1} \cos(\cdot) + \eqref{eq:zMode_newC_2} \sin(\cdot)$ and $\eqref{eq:zMode_newC_1} \sin(\cdot) - \eqref{eq:zMode_newC_2} \cos(\cdot)$, we find the differential equations for $a_1$ and $\phi$. Here, $(\cdot)$ denotes the phase
$t/\epsilon+\psi_1$.

\subsection{Recovering $a_0, a_1$ and $\phi_1$ from data}
\label{app:recover_0eig_ampl_phases}

We have computed a short time simulation at the grid points $t_i$ which yields the data $y_{1,i}, y_{2,i}, y_{3,i}$. We solve a system of first-order ordinary differential equations with three components. The first index refers to the component of the solution, and the second index $i$ refers to the time grid point $t_i$. We can recover the solution in transformed coordinates at the grid points from the data exploiting the following relations:
\begin{align}
    a_{0,i} &= y_{3,i} \\
    a_{1,i} &= \frac{1}{\sqrt{2}} \sqrt{y_{1,i}^2 + y_{2,i}^2} \\
    \tan \left( \frac{t_i}{\epsilon} +
    \psi_{1,i} \right) &= \frac{y_{2,i}}{y_{1,i}} 
\end{align}

\subsection{Relation between $a_1(t)$ and $a_0(t)$}
\label{app:relation_a_a_0}

Assume that $a_0(t)> 0$ and $a_1(t)> 0$. From
\begin{equation}
    \frac{\dot a_1(t)}{a_1(t)} = \frac{\dot a_0(t)}{a_0(t)}, 
\end{equation}
it follows that
\begin{equation}
    \frac{\mathrm d}{\mathrm dt}\log(a_1(t)) = \frac{\mathrm d}{\mathrm dt}\log(a_0(t)) .
\end{equation}
Integration with respect to time yields
\begin{equation}
    \log(a_1(t)) = \log(a_0(t)) + C
\end{equation}
where $C$ is a time-independent constant. Exponentiating both sides, we obtain
\begin{equation}
    a_1(t) = a_0(t) e^C .
\end{equation}

\section{Details for van der Pol oscillator}
\label{app:van_der_Pol}

\subsection{Differential equations for amplitude and phase}
\label{subsec:Amplitude_Phase_Equs}

The van der Pol system is given by
\begin{equation}
    \vb{\dot y} = \frac{1}{\epsilon} \mathcal{L} \vb{y} + \begin{pmatrix}
        0 \\ g(\vb{y})
    \end{pmatrix} , \qquad g(\vb{y}) = (1-y_1^2)y_2 , \qquad
    \mathcal{L} = \begin{pmatrix}
        0 & 1 \\ -1 & 0
    \end{pmatrix}.
\end{equation}

We consider the following representation for the solution $\vb{y}(t)$
\begin{equation}
\label{eq:wkb_ansatz_app}
    \vb{y}(t) = a(t) e^{\frac{i}{\epsilon}\phi(t)} \vb{r} + a(t) e^{-\frac{i}{\epsilon}\phi(t)} \overline{\vb{r}} , \qquad \vb{r} = \frac{1}{\sqrt{2} }
     \begin{pmatrix} 1 \\ i \end{pmatrix}, \overline{\vb{r}} = \frac{1}{\sqrt{2} }
     \begin{pmatrix} 1 \\ -i \end{pmatrix} .
\end{equation}
with
\begin{equation}
    \phi (t) = t + \frac{1}{\epsilon} \psi (t)
\end{equation}
Computing the derivative of \eqref{eq:wkb_ansatz_app} gives
\begin{equation}
\label{eq:der_wkb_ansatz_app}
    \vb{\dot y} = \left( \dot a + \frac{i}{\epsilon} a \dot \phi \right) e^{\frac{i}{\epsilon}\phi(t)} \vb{r} +  \left( \dot a - \frac{i}{\epsilon} a \dot \phi \right) e^{-\frac{i}{\epsilon}\phi(t)} \overline{\vb{r}}
\end{equation}

Projection on the eigenvectors of the linear operator yields
\begin{equation}
\label{eq:der_saclar_r}
    \vb{\dot y} \cdot \overline{\vb{r}} =   \left( \dot a + \frac{i}{\epsilon} a \dot \phi \right) e^{\frac{i}{\epsilon}\phi(t)} ; \qquad \vb{\dot y} \cdot \vb{r} = \left( \dot a - \frac{i}{\epsilon} a \dot \phi \right) e^{-\frac{i}{\epsilon}\phi(t)} .
\end{equation}

Using the representation \eqref{eq:wkb_ansatz_app} for $\vb{y}$ and applying it to the linear operator, we find
\begin{equation}
    \frac{1}{\epsilon} \mathcal{L} \vb{y} = \frac{i}{\epsilon} a e^{\frac{i}{\epsilon}\phi} \vb{r} - \frac{i}{\epsilon} a e^{-\frac{i}{\epsilon}\phi} \overline{\vb{r}} .
\end{equation}

Again we apply projection on the eigenvectors
\begin{equation}
\label{eq:L_scalar_r}
     \frac{1}{\epsilon} \mathcal{L} \vb{y} \cdot \overline{\vb{r}} = \frac{i}{\epsilon} a e^{\frac{i}{\epsilon}\phi}, \qquad
      \frac{1}{\epsilon} \mathcal{L} \vb{y} \cdot \vb{r} = -\frac{i}{\epsilon} a e^{-\frac{i}{\epsilon}\phi} .
\end{equation}

We use the identity
\begin{equation}
    \vb{\dot y} \cdot \overline{\vb{r}} = \frac{1}{\epsilon} \mathcal{L} \vb{y} \cdot \overline{\vb{r}} + \begin{pmatrix}
        0 \\ g(\vb{y})
    \end{pmatrix} \cdot \overline{\vb{r}}
\end{equation}

and plug in the relations \eqref{eq:der_saclar_r} and \eqref{eq:L_scalar_r} 
\begin{equation}
    \left( \dot a + \frac{i}{\epsilon} a \dot \phi \right) e^{\frac{i}{\epsilon}\phi} = \frac{i}{\epsilon} a e^{\frac{i}{\epsilon}\phi} + \begin{pmatrix}
        0 \\ g(\vb{y})
    \end{pmatrix} \cdot \overline{\vb{r}} .
\end{equation}

The computations are repeated with the other eigenvector
\begin{equation}
    \vb{\dot y} \cdot \vb{r} = \frac{1}{\epsilon} \mathcal{L} \vb{y} \cdot \vb{r} + \begin{pmatrix}
        0 \\ g(\vb{y})
    \end{pmatrix} \cdot \vb{r}
\end{equation}

and \eqref{eq:der_saclar_r} and \eqref{eq:L_scalar_r} are used again
\begin{equation}
    \left( \dot a - \frac{i}{\epsilon} a \dot \phi \right) e^{-\frac{i}{\epsilon}\phi} = -\frac{i}{\epsilon} a e^{-\frac{i}{\epsilon}\phi} + \begin{pmatrix}
        0 \\ g(\vb{y})
    \end{pmatrix} \cdot \vb{r} .
\end{equation}

Thus, we obtain the following system of equations for $\dot a$ and $\dot \phi$
\begin{align}
    \dot a + \frac{i}{\epsilon} a \dot \phi &= \frac{i}{\epsilon} a +  e^{-\frac{i}{\epsilon}\phi} \begin{pmatrix}
        0 \\ g(\vb{y})
    \end{pmatrix} \cdot \overline{\vb{r}}\\
    \dot a - \frac{i}{\epsilon} a \dot \phi &= -\frac{i}{\epsilon} a +  e^{\frac{i}{\epsilon}\phi} \begin{pmatrix}
        0 \\ g(\vb{y})
    \end{pmatrix} \cdot \vb{r} .
\end{align}

Algebraic manipulations lead to
\begin{align}
    2 \dot a &=   \begin{pmatrix}
        0 \\ g(\vb{y})
    \end{pmatrix} \left(  e^{-\frac{i}{\epsilon}\phi}  \overline{\vb{r}} +  e^{\frac{i}{\epsilon}\phi}  \vb{r} \right)\\
    2 \frac{i}{\epsilon} a \dot \phi &= 2 \frac{i}{\epsilon} a + \begin{pmatrix}
        0 \\ g(\vb{y})
    \end{pmatrix} \left(  e^{-\frac{i}{\epsilon}\phi}  \overline{\vb{r}} -  e^{\frac{i}{\epsilon}\phi}  \vb{r} \right) .
\end{align}

\subsection{Details of the right-hand side}
\label{subsec:furhter_details_rhs}

The function $g(\vb{y})$, which defines the nonlinearity of the van der Pol system, is given by
\begin{equation}
    g(\vb{y}) = (1-y_1^2) y_2 = y_2 - y_1^2 y_2 .
\end{equation}

The nonlinearity shall be expressed with the representation \eqref{eq:wkb_ansatz_app}. In the first step, we compute
\begin{align}
    y_1^2 &= \left( \frac{1}{\sqrt{2}} a e^{\frac{i}{\epsilon}\phi} + \frac{1}{\sqrt{2}} a e^{-\frac{i}{\epsilon}\phi} \right)^2 \\
    &= \frac{1}{2} a^2 \left(  e^{\frac{i}{\epsilon}\phi} + e^{-\frac{i}{\epsilon}\phi} \right)^2 \\
    &= \frac{1}{2} a^2 \left(  e^{2\frac{i}{\epsilon}\phi} + 2 + e^{-2\frac{i}{\epsilon}\phi} \right) .
\end{align}

In the second step, we make the following computation
\begin{align}
    y_1^2 y_2 &= \frac{1}{2} a^2 \left(  e^{2\frac{i}{\epsilon}\phi} + 2 + e^{-2\frac{i}{\epsilon}\phi} \right)^2 \frac{1}{\sqrt{2}} a i \left( e^{\frac{i}{\epsilon}\phi} - e^{-\frac{i}{\epsilon}\phi} \right)\\
    &= \frac{1}{2\sqrt{2}} i a^3 \left( e^{3\frac{i}{\epsilon}\phi} - e^{\frac{i}{\epsilon}\phi} + 2e^{\frac{i}{\epsilon}\phi} - 2 e^{-\frac{i}{\epsilon}\phi} + e^{-\frac{i}{\epsilon}\phi} - e^{-3\frac{i}{\epsilon}\phi}\right) \\
    &= \frac{1}{2\sqrt{2}} i a^3 \left( e^{3\frac{i}{\epsilon}\phi} + e^{\frac{i}{\epsilon}\phi} - e^{-\frac{i}{\epsilon}\phi} - e^{-3\frac{i}{\epsilon}\phi} \right) .
\end{align}

Now, the the projection on $\vb{r}$ is computed
\begin{align}
    \begin{pmatrix}
        0 \\ g(\vb{y}) 
    \end{pmatrix} \cdot \overline{\vb{r}} &= 
    \begin{pmatrix}
        0 \\ g(\vb{y}) 
    \end{pmatrix} \cdot \frac{1}{\sqrt{2}} \begin{pmatrix}
        1 \\ -i
    \end{pmatrix} =
    -\frac{i}{\sqrt{2}} g(\vb{y}) \\
    &= -\frac{i}{\sqrt{2}} \left( \frac{ia}{\sqrt{2}} \left( e^{\frac{i}{\epsilon}\phi} - e^{-\frac{i}{\epsilon}\phi} \right) -  \frac{1}{2\sqrt{2}} i a^3 \left( e^{3\frac{i}{\epsilon}\phi} + e^{\frac{i}{\epsilon}\phi} - e^{-\frac{i}{\epsilon}\phi} - e^{-3\frac{i}{\epsilon}\phi} \right)\right) \\
    &= \frac{1}{2} a \left( e^{\frac{i}{\epsilon}\phi} - e^{-\frac{i}{\epsilon}\phi} \right) - \frac{1}{4} a^3 \left( e^{3\frac{i}{\epsilon}\phi} + e^{\frac{i}{\epsilon}\phi} - e^{-\frac{i}{\epsilon}\phi} - e^{-3\frac{i}{\epsilon}\phi} \right),
\end{align}

followed by the projection on $\overline{\vb{r}}$
\begin{align}
 \begin{pmatrix}
        0 \\ g(\vb{y}) 
    \end{pmatrix} \cdot \vb{r} &= 
    \begin{pmatrix}
        0 \\ g(\vb{y}) 
    \end{pmatrix} \cdot \frac{1}{\sqrt{2}} \begin{pmatrix}
        1 \\ i
    \end{pmatrix} =
    \frac{i}{\sqrt{2}} g(\vb{y}) \\
    &= \frac{i}{\sqrt{2}} \left( \frac{ia}{\sqrt{2}} \left( e^{\frac{i}{\epsilon}\phi} - e^{-\frac{i}{\epsilon}\phi} \right) -  \frac{1}{2\sqrt{2}} i a^3 \left( e^{3\frac{i}{\epsilon}\phi} + e^{\frac{i}{\epsilon}\phi} - e^{-\frac{i}{\epsilon}\phi} - e^{-3\frac{i}{\epsilon}\phi} \right)\right) \\
    &= -\frac{1}{2} a \left( e^{\frac{i}{\epsilon}\phi} - e^{-\frac{i}{\epsilon}\phi} \right) + \frac{1}{4} a^3 \left( e^{3\frac{i}{\epsilon}\phi} + e^{\frac{i}{\epsilon}\phi} - e^{-\frac{i}{\epsilon}\phi} - e^{-3\frac{i}{\epsilon}\phi} \right) .
\end{align}

For the right-hand side of $\dot a$, we get
\begin{align}
\begin{split}
  \frac{1}{2} \begin{pmatrix}
        0 \\ g(\vb{y}) 
    \end{pmatrix} \cdot ( e^{\frac{i}{\epsilon}\phi} \vb{r} +  e^{-\frac{i}{\epsilon}\phi} \overline{\vb{r}}) 
        &= \frac{1}{2} \left( - \frac{1}{2} a \left( e^{2\frac{i}{\epsilon}\phi} -1 \right) + \frac{1}{4} a^3 \left(  e^{4\frac{i}{\epsilon}\phi} +  e^{2\frac{i}{\epsilon}\phi} -1 -  e^{-2\frac{i}{\epsilon}\phi} \right) \right) \\
        & \ \quad + \frac{1}{2} \left( \frac{1}{2} a (1- e^{-2\frac{i}{\epsilon}\phi}) -\frac{1}{4} a^3 \left(  e^{2\frac{i}{\epsilon}\phi} +1 -  e^{-2\frac{i}{\epsilon}\phi} -  e^{-4\frac{i}{\epsilon}\phi} \right) \right) \\
        &= \frac{1}{4} a - \frac{1}{8} a^3 + \frac{1}{4} a - \frac{1}{8} a^3 \\
        & \ \quad  -\frac{1}{4} a e^{2\frac{i}{\epsilon}\phi} + \frac{1}{8} a^3 e^{4\frac{i}{\epsilon}\phi} + \frac{1}{8} a^3 e^{2\frac{i}{\epsilon}\phi} - \frac{1}{8} a^3 e^{-2\frac{i}{\epsilon}\phi} \\
        & \ \quad -\frac{1}{4} a e^{- 2\frac{i}{\epsilon}\phi} -\frac{1}{8} a^3 e^{ 2\frac{i}{\epsilon}\phi} + \frac{1}{8} a^3 e^{- 2\frac{i}{\epsilon}\phi} + \frac{1}{8} a^3 e^{- 4\frac{i} {\epsilon}\phi} \\
        &= \frac{1}{2} a - \frac{1}{4} a^3 - \frac{1}{4} a e^{2\frac{i}{\epsilon}\phi} - \frac{1}{4} a e^{- 2\frac{i}{\epsilon}\phi} + \frac{1}{8} a^3 e^{4\frac{i}{\epsilon}\phi} + \frac{1}{8} a^3 e^{- 4\frac{i}{\epsilon}\phi} .
    \end{split} 
\end{align}

The right-hand side of $\dot \phi$ contains the following terms
\begin{equation}
    \begin{split}
        \begin{pmatrix}
        0 \\ g(\vb{y}) 
    \end{pmatrix} \cdot ( e^{-\frac{i}{\epsilon}\phi} \overline{\vb{r}} -  e^{\frac{i}{\epsilon}\phi} \vb{r}) 
        &=  \frac{1}{2} a \left(1 - e^{-2\frac{i}{\epsilon}\phi}  \right) - \frac{1}{4} a^3 \left(  e^{2\frac{i}{\epsilon}\phi} +1 -  e^{-2\frac{i}{\epsilon}\phi} -  e^{-4\frac{i}{\epsilon}\phi} \right)  \\
        & \ \quad -
        \left( - \frac{1}{2} a \left( e^{2\frac{i}{\epsilon}\phi} -1 \right) + \frac{1}{4} a^3 \left(  e^{4\frac{i}{\epsilon}\phi} +  e^{2\frac{i}{\epsilon}\phi} - 1 - e^{-2\frac{i}{\epsilon}\phi} \right)\right) \\ 
        &= \frac{1}{2} a - \frac{1}{2} a - \frac{1}{4}a ^3 + \frac{1}{4} a^3 - \frac{1}{2} a e^{-2\frac{i}{\epsilon}\phi} +  \frac{1}{2} a e^{2\frac{i}{\epsilon}\phi}\\
        & \ \quad - \frac{1}{4} a^3 e^{2\frac{i}{\epsilon}\phi} + \frac{1}{4} a^3 e^{-2\frac{i}{\epsilon}\phi} + \frac{1}{4} a^3 e^{-4\frac{i}{\epsilon}\phi}  \\
         & \ \quad -
        \frac{1}{4} a^3 e^{4\frac{i}{\epsilon}\phi} - \frac{1}{4} a^3 e^{2\frac{i}{\epsilon}\phi} + \frac{1}{4} a^3 e^{-2\frac{i}{\epsilon}\phi} \\
        &= - \frac{1}{2} a e^{-2\frac{i}{\epsilon}\phi} +  \frac{1}{2} a e^{2\frac{i}{\epsilon}\phi} -  \frac{1}{4} a^3 e^{4\frac{i}{\epsilon}\phi} +  \frac{1}{4} a^3 e^{-4\frac{i}{\epsilon}\phi} - \frac{1}{2} a^3 e^{2\frac{i}{\epsilon}\phi} + \frac{1}{2} a^3 e^{-2\frac{i}{\epsilon}\phi} .
    \end{split}
\end{equation}

\subsection{Asymptotics for van der Pol}
\label{app:Asymptotics_for_Van_der_Pol}

Using the results of \cref{subsec:Amplitude_Phase_Equs} and \cref{subsec:furhter_details_rhs}, we find
\begin{align}
    \dot a &= \frac{1}{2} a - \frac{1}{4} a^3 - \frac{1}{2} a \cos(\frac{2}{\epsilon}(t+\epsilon\psi)) + \frac{1}{4} a^3 \cos(\frac{4}{\epsilon} (t+\epsilon\psi)) \\
    \dot \psi &= \frac{1}{2} \sin(\frac{2}{\epsilon}(t+\epsilon\psi)) - \frac{1}{2} a^2 \sin(\frac{2}{\epsilon}(t+\epsilon\psi)) - \frac{1}{4} a^2 \sin(\frac{4}{\epsilon}(t+\epsilon\psi))
\end{align}

We assume that $a$ and $\psi$ depend on two time scales
\begin{equation}
    a(t) = a\left(t, \frac{t}{\epsilon} \right) = a(t,\tau) , \qquad \qquad \psi(t) = \psi\left(t, \frac{t}{\epsilon} \right) = \psi(t, \tau) .
\end{equation}

We choose the following asymptotic expansions:
\begin{align}
    a(t) &= a^{(0)}(t) + \epsilon a^{(1)}(t, \tau) + \epsilon^2 a^{(2)}(t, \tau) + o(\epsilon^2) \\
    \psi(t) &= \psi^{(0)}(t) + \epsilon \psi^{(1)}(t, \tau) + \epsilon^2 \psi^{(2)}(t, \tau) + o(\epsilon^2) 
\end{align}

For the asymptotic expansion of $\psi$ we find:
\begin{align}
   & \psi^{(0)}_t + \epsilon \psi^{(1)}_t + \psi^{(1)}_{\tau} + \epsilon \psi^{(2)}_{\tau} + o(\epsilon) \\
  &  = \frac{1}{2} \sin(2\tau + 2\psi^{(0)} + 2\epsilon\psi^{(1)}+ o(\epsilon)) \\
  & - \frac{1}{2} ({a^{(0)}}^2 + 2 \epsilon a^{(0)} a^{(1)} + o(\epsilon))   \sin(2\tau + 2\psi^{(0)} + 2\epsilon\psi^{(1)}+ o(\epsilon)) \\
  &  - \frac{1}{4}  (a{^{(0)}}^2 + 2\epsilon a^{(0)} a^{(1)} + o(\epsilon)) \sin(4\tau + 4\psi^{(0)} + 4\epsilon\psi^{(1)}+ o(\epsilon)) \\
  & \approx \Bigg( \frac{1}{2} \sin(2\tau) \left( \cos(2 \psi^{(0)}) - 2 \epsilon \psi^{(1)} \sin(2\psi^{(0)}) \right) \\
  & + \frac{1}{2} \cos(2\tau) \left( \sin(2 \psi^{(0)})  + 2 \epsilon \psi^{(1)} \sin(2\psi^{(0)}) \right) \Bigg) \\
  & \times \left( 1 - {a^{(0)}}^2 - 2 \epsilon a^{(0)} a^{(1)} + o(\epsilon) \right) \\
  &- \frac{1}{4} \Bigg(\cos(4\tau) \left( \sin(4 \psi^{(0)})  + 4 \epsilon \psi^{(1)} \cos(4\psi^{(0)}) \right) \\
  & \sin(4\tau) \left( \cos(4 \psi^{(0)})  - 4 \epsilon \psi^{(1)} \sin(4\psi^{(0)}) \right) \Bigg) \\
  & \times \left( {a^{(0)}}^2 + 2 \epsilon a^{(0)} a^{(1)} + o(\epsilon) \right) 
\end{align}

We used
\begin{align}
    \sin(2\tau + 2\psi^{(0)} + 2\epsilon \psi^{(1)}+ o(\epsilon)) &= \sin(2\tau) \cos( 2\psi^{(0)} + 2\epsilon\psi^{(1)}+ o(\epsilon) ) \\
    & \quad + \cos(2\tau) \sin( 2\psi^{(0)} + 2\epsilon\psi^{(1)}+ o(\epsilon) ) \\ 
     & \quad \approx \sin(2\tau) \left( \cos(2\psi^{(0)}) - 2 \epsilon \psi^{(1)} \sin(2\psi^{(0)}) \right) \\
     & \quad + \cos(2\tau) \left( \sin(2\psi^{(0)}) + 2 \epsilon \psi^{(1)} \cos(2\psi^{(0)}) \right)
\end{align}

The $O(1)$ equation of $\psi$ is given by
\begin{align}
    \psi_t^{(0)} + \psi_{\tau}^{(1)} & =
    \left(1 - {a^{(0)}}^2 \right) 
    \left(\frac{1}{2} \sin(2\tau)  \cos(2 \psi^{(0)})
    +\frac{1}{2} \cos(2\tau)  \sin(2 \psi^{(0)})\right) \\
    & \quad - \frac{1}{4} {a^{(0)}}^2
     \left( \sin(4\tau)  \cos(4 \psi^{(0)})
    + \cos(4\tau)  \sin(4 \psi^{(0)})\right)
\end{align}

Integrating with respect to $\tau$ and applying the sublinear growth condition, see \textit{e.g} \cite{Klein2010MultipleScales}, yields\footnote{Integrating the equation with respect to $\tau$ yields the term $\tau \psi^{(0)}_{t}(t)$. To avoid secular growth, we require $\psi^{(0)}_{t}(t) = 0$.}
\begin{equation}
    \psi_t^{(0)} = 0
\end{equation}


For the amplitude equation we find:
\begin{align}
    & a^{(0)}_t + \epsilon a^{(1)}_t +a^{(1)}_{\tau} + \epsilon a^{(2)}_{\tau} + o(\epsilon) \\
    & \quad = \frac{1}{2} \left( a^{(0)} + \epsilon  a^{(1)} + o(\epsilon) \right) - \frac{1}{4} \left(  {a^{(0)}}^3 + 3 \epsilon { a^{(0)}}^2
    a^{(1)} + o(\epsilon) \right) \\
    & \quad - \frac{1}{2} \left( a^{(0)} + \epsilon  a^{(1)} + o(\epsilon) \right) \cos(2 \tau + 2 \psi^{(0)} + 2 \epsilon \psi^{(1)} + o(\epsilon)) \\
    & \quad + \frac{1}{4} \left(  {a^{(0)}}^3 + 3 \epsilon { a^{(0)}}^2
    a^{(0)} + o(\epsilon) \right) \cos(4 \tau + 4 \psi^{(0)} + 4 \epsilon \psi^{(1)} + o(\epsilon)) \\
    & \quad \approx \frac{1}{2} \left( a^{(0)} + \epsilon  a^{(1)} + o(\epsilon) \right) - \frac{1}{4} \left(  {a^{(0)}}^3 + 3 \epsilon { a^{(0)}}^2
    a^{(1)} + o(\epsilon) \right) \\
    & \quad - \frac{1}{2} \left( a^{(0)} + \epsilon  a^{(1)} + o(\epsilon) \right) \\
    & \quad \times \left(  \cos(2\tau) \left( \cos( 2 \psi^{(0)}) - 2 \epsilon \psi^{(1)} \sin(2 \psi^{(0)})) \right) 
    - \sin(2\tau) \left( \sin( 2 \psi^{(0)}) + 2 \epsilon \psi^{(1)} \cos(2 \psi^{(0)})) \right) \right) \\
    & \quad + \frac{1}{4} \left(  {a^{(0)}}^3 + 3 \epsilon { a^{(0)}}^2  a^{(0)} + o(\epsilon) \right) \\
    & \quad \times \left(  \cos(4\tau) \left( \cos(4 \psi^{(0)}) - 4 \epsilon \psi^{(1)} \sin(4 \psi^{(0)})) \right) 
    - \sin(4\tau) \left( \sin(4 \psi^{(0)}) + 4 \epsilon \psi^{(1)} \cos(4 \psi^{(0)})) \right) \right)
\end{align}

We used
\begin{equation}
    \left( a^{(0)} + \epsilon  a^{(1)} + o(\epsilon) \right)^3 = {a^{(0)}}^3 + 3 \epsilon { a^{(0)}}^2
    a^{(0)} + o(\epsilon) .
\end{equation}

and
\begin{align}
   & \cos(2 \tau + 2 \psi^{(0)} + 2 \epsilon \psi^{(1)})  \\
    & = \cos(2\tau) \left( \cos( 2 \psi^{(0)} + 2 \epsilon \psi^{(1)}) \right) - \sin(2\tau) \left( \sin( 2 \psi^{(0)} + 2 \epsilon \psi^{(1)}) \right) \\
    & \approx \cos(2\tau) \left( \cos( 2 \psi^{(0)}) - 2 \epsilon \psi^{(1)} \sin(2 \psi^{(0)})) \right) 
    - \sin(2\tau) \left( \sin( 2 \psi^{(0)}) + 2 \epsilon \psi^{(1)} \cos(2 \psi^{(0)})) \right) 
\end{align}

The $O(1)$ equation is given by:
\begin{align}
    a^{(0)}_t + a^{(1)}_{\tau} &= \frac{1}{2} a^{(0)} - \frac{1}{4} {a^{(0)}}^3 - \frac{1}{2} a^{(0)} \left( \cos(2\tau) \cos(2 \psi^{(0)}) - \sin(2\tau) \sin(2 \psi^{(0)}) \right) \\
    & + \frac{1}{4} {a^{(0)}}^3 \left( \cos(4\tau) \cos(4 \psi^{(0)}) - \sin(4\tau) \sin(4 \psi^{(0)}) \right)
\end{align}


We integrate the $O(1)$ equation with respect to $\tau$. To avoid secular growth (sublinear growth condition), we require
\begin{equation}
a_t^{(0)} = \frac{1}{2} a^{(0)} - \frac{1}{4} \big(a^{(0)}\big)^3,
\end{equation}
which determines the leading-order asymptotic behavior of $a(t)$.


\section{Pseudocode}
\label{app:pseudocode}

\begin{algorithm}[H]
\caption{RHS of van der Pol System}
\begin{algorithmic}[1]
\Function{RHS\_vdP}{$x, \varepsilon$}
    \State $x_1 \gets x[0]$
    \State $x_2 \gets x[1]$
    \State $\dot{x}_1 \gets x_2 / \varepsilon$
    \State $\dot{x}_2 \gets -x_1 / \varepsilon + (1 - x_1^2) \cdot x_2$
    \State \Return $[\dot{x}_1, \dot{x}_2]$
\EndFunction
\end{algorithmic}
\end{algorithm}

\begin{algorithm}[H]
\caption{RHS for system with zero eigenvalue}
\begin{algorithmic}[1]
\Function{LinearOperator}{$x, \varepsilon$}
    \State $(y_1, y_2, y_3) \gets x$
    \State \Return $\frac{1}{\varepsilon} \cdot [\cdot y_2, - y_1, 0]$
\EndFunction

\Function{NonlinearTerm}{$x, \lambda$}
    \State $(y_1, y_2, y_3) \gets x$
    \State $sq \gets y_1^2 + y_2^2 + y_3^2$
    \State \Return $(-\lambda) \cdot sq \cdot [y_1, y_2, y_3]$
\EndFunction

\Function{RHS}{$x, \varepsilon, \lambda$}
    \State $L \gets$ LinearOperator$(x, \varepsilon, B)$
    \State $N \gets$ NonlinearTerm$(x, \lambda)$
    \State \Return $L + N$  \Comment{Sum of linear and nonlinear terms}
\EndFunction
\end{algorithmic}
\end{algorithm}
 
\begin{algorithm}[H]
\caption{Fourth-order Runge--Kutta time-stepping scheme (RK4)}
\begin{algorithmic}[1]
\Function{RK4Step}{$\text{rhs}, x, \Delta t$}
  \State $k_1 \gets \text{rhs}(x)$
  \State $k_2 \gets \text{rhs}(x + 0.5 \cdot \Delta t \cdot k_1)$
  \State $k_3 \gets \text{rhs}(x + 0.5 \cdot \Delta t \cdot k_2)$
  \State $k_4 \gets \text{rhs}(x + \Delta t \cdot k_3)$
  \State \Return $x + \frac{\Delta t}{6} (k_1 + 2k_2 + 2k_3 + k_4)$
\EndFunction
\end{algorithmic}
\end{algorithm}

\begin{algorithm}[H]
\caption{Time integration using RK4}
\begin{algorithmic}[1]
\Function{Integrate}{$\text{rhs}, x_0, t_{\min}, t_{\max}, \Delta t$}
  \State $t \gets$ array from $t_{\min}$ to $t_{\max}$ with spacing $\Delta t$
  \State $L \gets$ length of $t$
  \State $\text{dim} \gets$ length of $x_0$
  \State $x \gets$ zero array of shape $(\text{dim}, L)$
  \State $x[:,0] \gets x_0$
  \For{$n = 0$ to $L - 2$}
    \State $x[:,n+1] \gets$ \Call{RK4Step}{$\text{rhs}, x[:,n], \Delta t$}
  \EndFor
  \State \Return $t, x$
\EndFunction
\end{algorithmic}
\end{algorithm}

\begin{algorithm}[H]
\caption{Find Amplitude}
\begin{algorithmic}[1]
\Function{FindAmplitude}{$y_1^{data}, y_2^{data}$}
    \State $amplitude \gets \frac{1}{\sqrt{2}} \cdot \sqrt{(y_1^{data})^2 + (y_2^{data})^2}$
    \State \textbf{Return:} $amplitude$
\EndFunction
\end{algorithmic}
\end{algorithm}

\begin{algorithm}[H]
\caption{Compute phase $\phi$ from data}
\begin{algorithmic}[1]
\Function{FindPhi}{$y_1^\text{data}, y_2^\text{data}, \varepsilon$}
    \State \(\triangleright\) Compute wrapped phase in $(-\pi, \pi]$ from Cartesian data
    \State $\phi_{\varepsilon}^{\text{wrapped}} \gets \arctan2(y_2^\text{data}, y_1^\text{data})$
    \State \(\triangleright\) Remove $2\pi$ discontinuities to get continuous phase
    \State $\phi \gets \varepsilon \cdot \operatorname{Unwrap}(\phi_{\varepsilon}^{\text{wrapped}})$
    \State \textbf{return} $\phi$
\EndFunction
\end{algorithmic}
\end{algorithm}




\begin{algorithm}[H]
\caption{Extrapolate van der Pol system}
\begin{algorithmic}[1]
\Function{Extrapolate}{$a\_, b\_, c\_, d\_, bridge\_grid\_, \varepsilon$}
    \State $appr\_a \gets a\_ + b\_ \cdot bridge\_grid\_$
    \State $appr\_\phi \gets c\_ + d\_ \cdot bridge\_grid\_$
    \State $y1\_ext\_ \gets \sqrt{2} \cdot appr\_a \cdot \cos(appr\_\phi / \varepsilon)$
    \State $y2\_ext\_ \gets (-\sqrt{2}) \cdot appr\_a \cdot \sin(appr\_\phi / \varepsilon)$
    \State \textbf{Return:} $y1\_ext\_,$ $y2\_ext\_$
\EndFunction
\end{algorithmic}
\end{algorithm}

\begin{algorithm}[H]
\caption{Scaled Kernel Function}
\begin{algorithmic}[1]
\Function{scaled\_rho}{$\eta, dt$}
    \State $\rho_0 \gets 0.007029858406609657$ \Comment{ constant normalizes kernel in $l_1$ norm}
    \State $K \gets \eta / dt$ \Comment{number of smoothing points}
    \State $grid_i \gets -0.5 \cdot \eta + i \cdot dt$, \quad $i = 0, \dots, K$ \Comment{grid in $[-\eta/2, \eta/2]$}
    \State $s \gets grid / \eta$  \Comment{grid in $[-1/2, 1/2]$}
    \State $\rho \gets$ zeros\_like$(s)$
    \For{each $i$ where $-0.5 < |s[i]| < 0.5$}
        \State $\rho[i] \gets 1/\rho_0 \cdot  \exp\left(1 / ((s[i]-0.5) \cdot (s[i]+0.5))\right)$ \Comment{evaluate $\rho(s_i)$}
    \EndFor
    \State $sc\_\rho \gets \rho / \eta$ \Comment{evaluate $1/\eta \cdot \rho(s_i/\eta)$}
    \State \Return $grid, sc\_\rho$
\EndFunction
\end{algorithmic}
\end{algorithm}

\bibliographystyle{plain} 
\bibliography{references} 


\end{document}